\documentclass[journal ]{new-aiaa}
\usepackage[utf8]{inputenc}
\usepackage{textcomp}

\usepackage{graphicx}
\usepackage{amsmath}
\usepackage[version=4]{mhchem}
\usepackage{siunitx}
\usepackage{longtable,tabularx}
\usepackage{subcaption}

\usepackage{multirow, makecell, diagbox, float} 

\newtheorem{theorem}{Theorem} 

\title{Attitude Motion of Unbalanced Partial-Spin Spacecraft}

\author{Jingyuan Wu \footnote{Master student, Institute of Mechanics, wujingyuan@imech.ac.cn.}, Wenhao Li \footnote{Professor, Institute of Mechanics, liwenhao@imech.ac.cn.} and Guanhua Feng \footnote{Assistant professor, Institute of Mechanics, fengguanhua@imech.ac.cn.}}
\affil{ Chinese Academy of Sciences, Beijing, 100190}

\begin{document}

\maketitle

 \begin{abstract}
  This paper focuses on the attitude motion of spacecraft systems featuring asymmetric spacecraft platforms and unbalanced rotors. Through perturbation expansion, the spacecraft dynamic equation is simplified as a linear periodically time-varying (LPTV) system, and analytical solutions of angular velocities are derived. A novel stability criterion is proposed, providing insights into stability transitions and dynamic behavior. Furthermore, the characteristic motions of precession and nutation are analytically investigated, highlighting their dependence on system inertia properties. Numerical simulations are performed to validate the analytical results, showing excellent agreement within defined operational domains and offering error bounds for linearization. The findings advance the understanding of unbalanced partial-spin spacecraft dynamics and provide practical guidance for the design and optimization of such systems.
 
\end{abstract}

\section*{Nomenclature}


{\renewcommand\arraystretch{1.0}
\noindent\begin{longtable*}{@{}l @{\quad=\quad} l@{}}
$\boldsymbol{H}_{\Sigma }^S$ & angular momentum of the spacecraft system with respect to its mass center, \si{kg \cdot m^2\cdot s^{-1}}  \\
$\boldsymbol{I}_{A}$  &  inertia tensor of body A with respect to its mass center, \si{kg\cdot m^2}\\
$\boldsymbol{I}_{B}$  &  inertia tensor of body B with respect to its mass center, \si{kg\cdot m^2}\\
$I_{B_Y}$, $I_{B_R}$ &  matrix elements of $I_{B}$, \si{kg\cdot m^2}\\
$\boldsymbol{I}_{S}$  &  inertia tensor of A and B's mass centers with respect to the mass center of the spacecraft, \si{kg\cdot m^2}\\
$\boldsymbol{I}_{1}$ & equivalent spin inertia tensor, \si{kg\cdot m^2}\\
$I_{xx}$, $I_{yy}$, $I_{zz}$, $I_{xy}$&matrix elements of $I_{1}$, \si{kg\cdot m^2}\\
$\boldsymbol{I}_{2}$ & equivalent spacecraft inertia tensor, \si{kg\cdot m^2}\\
$h$& distance between the A's mass center and the spin axis, m\\
$d$& distance between the B's mass center and the spin plane, m\\

$m_A$, $m_B$ & mass of A and B, kg\\

$ \boldsymbol{\omega}$  &   angular velocity of the spacecraft platform, rad/s\\
$ \boldsymbol{\Omega}$  &   relative angular velocity between the spacecraft platform and rotor, rad/s\\
$ \tau$  & dimensionless time\\
$\boldsymbol{\Omega}_0$  &$ \frac{\boldsymbol{\Omega} }{\left | \boldsymbol{\Omega}\right |  }$ , dimensionless relative angular velocity  between the spacecraft platform and rotor\\
$ \boldsymbol{T} $  & rotation matrix between the body-fixed frame of the spacecraft platform and the rotor\\
$\boldsymbol{\theta} $  & $(\theta_x,\theta_y,\theta_z)^T$, Euler angles of the spacecraft platform in Z-Y-X rotation sequence, rad\\

\end{longtable*}}

\section{Introduction}
 \lettrine{W}{ith} the increasing demand for space exploration and utilization, spin has expanded beyond traditional spacecraft spin stabilization to encompass diverse applications in specific space missions. These applications utilize rotating components (or rotors) to achieve intended objectives, some of which may cause the mass distribution changes, such as material transport in spinning systems\cite{feng2022theoretical}, stable deployment using deployable booms \cite{schenk2014review}, and directional observation with spinning camera payloads \cite{WANG201991, chandra2021review}, and more. Consequently, the stability of the system attitude and platform orientation is severely affected, even leading to mission failure. Therefore, a comprehensive investigation of spin spacecraft attitude dynamics considering rotor unbalance effects is of paramount importance.

    Extensive studies have been conducted on the attitude dynamics of spinning spacecraft. The foundational work dates back to 1958 when Bracewell's seminal study revealed that rigid satellites can only maintain stable rotation about the principal axis with the maximum moment of inertia in dissipative systems \cite{bracewell1958rotation}. This fundamental principle, known as the maximum axis theorem, was complemented by his introduction of the energy dissipation method to analyze mechanical energy's influence on the spin axis's orientation. This method offers significant advantages for long-term attitude analysis of various spacecraft configurations, including rigid \cite{spencer1974energy} and flexible bodies \cite{tsuchiya1974dynamics}, often providing intuitive geometric insights. However, the method inherently lacks mathematical rigor and struggles to quantify the rate of energy dissipation in the system accurately\cite{likins1966effects}.

    Dual-spin spacecraft, which come in various configurations, have also been the subject of studies on their attitude stability and analytical solutions. For dual-spin spacecraft, the common assumption is that the platform’s angular velocity is relatively small to ensure the precise orientation of onboard instruments, with flywheels or other rotors used for attitude stabilization. In such cases, the balance and axial symmetry of the rotor are critical for simplifying and solving the problem. The simplest scenario assumes the rotor is axially symmetric and the rotational axis passes through the mass center of each component, allowing analytical solutions to be directly derived. For instance, Cochran utilized Euler’s equations to derive analytical solutions for two Euler angles and the relative rotation angle, introducing the concept of a “trap state” where the nutation angle increases indefinitely under specific conditions \cite{cochran1982attitude}.

    For configurations involving damped springs, Likins applied Routh’s criteria and the energy-sink method to establish stability criteria \cite{likins2003attitude}. In cases with small external disturbances, Tao employed a time-scale approach to provide asymptotic solutions for satellite attitude angular velocities \cite{tao1975satellite}. For dual-spin spacecraft with asymmetric rotors, Kazuo Tsuchiya studied their dynamic behavior, analyzed stability using the averaging method, and developed criteria for entering a “trap state” \cite{tsuchiya1979attitude}. Similarly, M.P. Scher and colleagues investigated configurations where the rotor’s spin axis formed an angle with the platform’s principal axis, identifying two potential “trap states” for such asymmetrical dual-spin systems \cite{scher1974dynamic}. Many other studies have examined similar spacecraft configurations \cite{bainum1970motion, janssens2011stability, meng2014attitude, gasbarri2016dynamic,liu2018attitude}.

    For cases where the rotor is symmetric but unbalanced, general research methods typically model the system as a point mass attached to a symmetric rotor or introduce perturbation torques\cite{liu2013chaos}. However, when the rotor exhibits strong asymmetry, such as in satellites with spinning long rods or tethered satellite systems, it is crucial to consider the change of mass distribution and its impact on spacecraft attitude dynamics. In such cases, traditional approaches may no longer be applicable, and the stability of attitude motion and analytical solutions require further investigation.

    Despite extensive research on single-body and two-body spacecraft attitude computation and stability analysis, most existing studies impose requirements on spacecraft configurations. These analytical approaches rely on restrictive assumptions, such as balanced systems with coinciding centers of mass and spin axes, or axisymmetric geometric configurations. Consequently, deriving precise analytical solutions is still challenging for irregular and asymmetric spacecraft systems, which increasingly characterize modern space missions. Motivated by this, the dynamic behavior of a spacecraft featuring an asymmetric platform and an unbalanced asymmetric rotor is investigated in the paper. The contributions of this paper are summarized as follows:
    \begin{itemize}
        \item A dynamic model for asymmetric spacecraft with unbalanced rotors is proposed, incorporating both platform and rotational inertia coupling effects. The model is linearized into a linear periodically time-varying system and analytical solutions are derived.
        \item A stability criterion based on primary inertial properties is formulated for unbalanced partial-spin spacecraft, enabling precise predictions of stability transitions.
        \item The study provides a detailed analytical description of precession and nutation phenomena in asymmetric spacecraft systems, offering new insights into their dependence on system inertial properties.
    \end{itemize}
    
    The remainder of this paper is organized as follows. The nonlinear dynamic model of the spacecraft system is developed in Section \ref{sec2-modeling}. Then in Section \ref{sec3-analyticalSol-stability}, the model is simplified into a linear periodically time-varying (LPTV) system. Analytical solutions for the LPTV system are derived, along with the stability criteria for spacecraft attitude. In Section \ref{sec4-precession-nutation}, the precession and nutation motions of the spacecraft system are further analyzed. In Sections \ref{sec5-error-ayalysis} and \ref{sec6-examples}, the error analysis and numerical examples are conducted to assess the applicability and effectiveness of analytical solutions. Finally, the conclusions of the paper are listed in Section \ref{sec7-conclusion}.

\section{Dynamic Modeling}
\label{sec2-modeling}
    \subsection{Physical description of the spacecraft system}
    \label{subsec2-1}

    The partial-spin spacecraft configuration analyzed in this paper is shown in Fig. \ref{fig1-system-diagram}. The system comprises two parts: a non-axisymmetric platform and an unbalanced asymmetric rotor. The rotor spins around a fixed axis \(S_A\) with a relative angular velocity ${\boldsymbol{\Omega}}$. Since the rotation axis \(S_A\) does not pass through the mass center \(O_A\) of the rotor, the unbalanced motion of the rotor induces an angular velocity $\boldsymbol{\omega}$ in the spacecraft platform. Here, $ \boldsymbol{\omega}$ represents the angular velocity of the platform’s body-fixed coordinate system relative to the inertial coordinate system (the definitions are detailed in Section \ref{subsec2-2}).  

    The analysis assumes that the platform is initially stationary, i.e., $ {\boldsymbol{\omega}}|_{t=0} = 0$, and that the system evolves free from external torques during the subsequent motion. All components are modeled as rigid bodies.

    \subsection{{Definition of the coordinate system}}
    \label{subsec2-2}

\begin{figure}[h]
    \centering
    \includegraphics[width=0.9\textwidth]{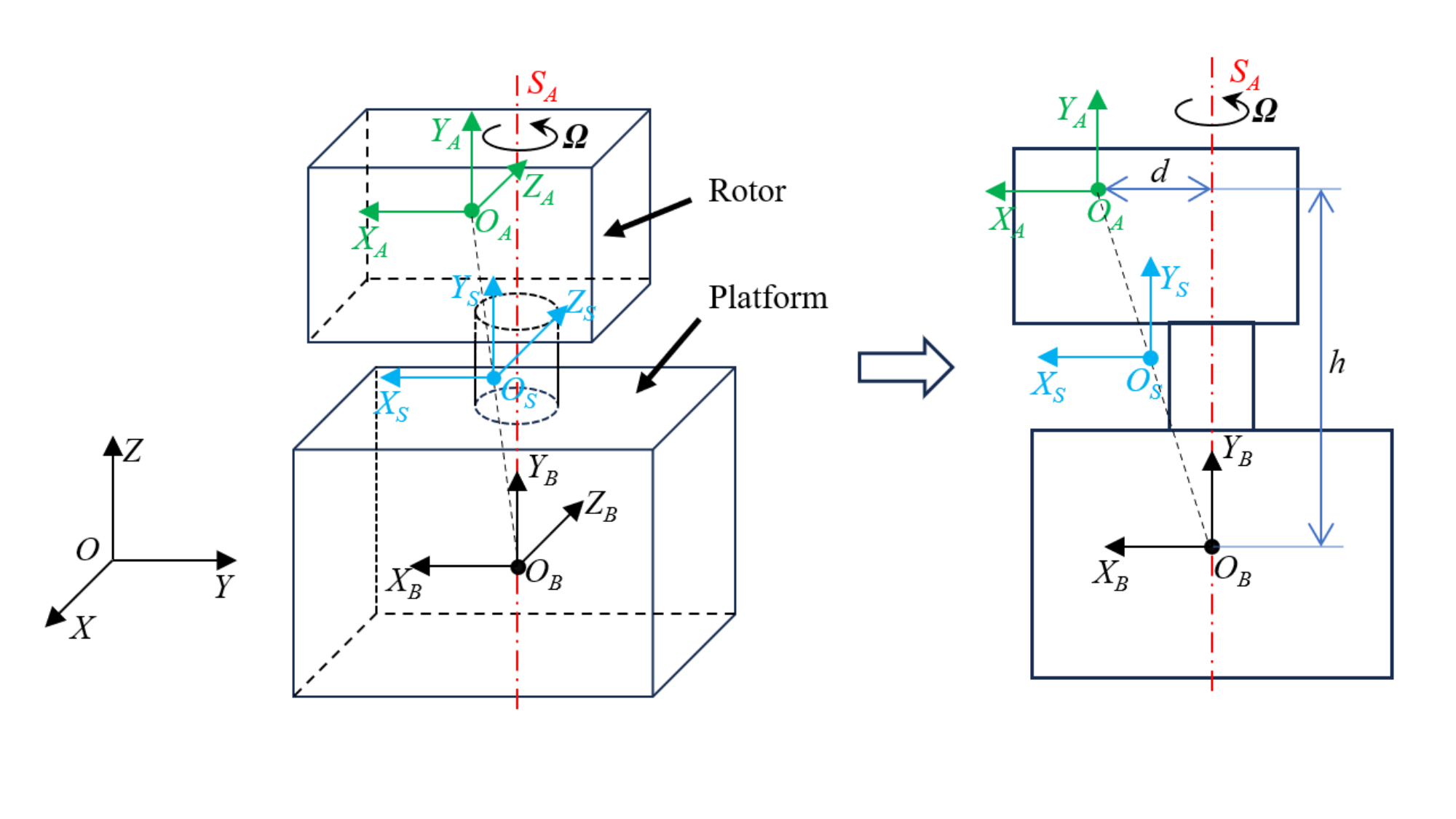}
    \caption{Partial spin spacecraft system}
    \label{fig1-system-diagram}
    \end{figure}

    The coordinate systems are defined herein, as depicted in Fig. \ref{fig1-system-diagram}.

    \begin{itemize}
     
        \item Body-fixed frame of the rotor \(O_A X_A Y_A Z_A\). The origin \(O_A\) is located at the mass center of the rotor. The three axes are aligned with the principal axes of inertia of the rotor and rotate synchronously with it. These axes adhere to the right-hand rule. Specifically, the \(Y_A\)-axis aligns counterclockwise with the rotation axis \(S_A\). The \(X_A\)-axis and \(Y_A\)-axis coplanar within the rotation plane of the rotor.
        \item Body-fixed frame of the platform \(O_B X_B Y_B Z_B\). The origin \(O_B\) is located at the mass center of the platform. The three axes are aligned with the principal axes of inertia of the platform and rotate together with it. These axes also satisfy the right-hand rule. Similarly, the \(Y_B\)-axis aligns counterclockwise with the rotation axis \(S_A\).
        \item Body-fixed frame of the entire system \(O_S X_S Y_S Z_S\). The origin \(O_S\) is located at the mass center of the entire system. The axes are oriented consistently with the platform's body-fixed coordinate system \(O_B X_B Y_B Z_B\). Initially, the \(X_A\)-axis, \(X_B\)-axis, and \(X_S\)-axis are collinear.
        \item Inertial frame \(O_SXYZ\). The origin \(O\) is located at the mass center of the entire system. The \(Y\)-axis, along with the \(X\)-axis and \(Z\)-axis, forms a right-handed coordinate system. At the beginning, the $X$, $Y$, and $Z$ axes are aligned with the $X_S$, $Y_S$, and $Z_S$ axes.
    \end{itemize}

    \subsection{{Dynamic modeling of spacecraft system}}
    \label{subsec2-3}
           Traditionally, unbalanced spacecraft systems have been simplified and modeled as balanced systems with discrete mass points. However, this approximation fails to capture critical coupling effects arising from mass distribution asymmetries—effects that significantly impact system dynamics in practical applications. This paper presents an enhanced model that treats the spacecraft rotor as an asymmetric rigid body with its center of mass displaced from its geometric spin axes, thereby incorporating these crucial dynamic interactions.
           According to König’s theorem for angular momentum, the total angular momentum of the system can be decomposed into four parts:
           
        (1) The angular momentum of the platform with respect to its own center of mass:
        \begin{equation}
        \boldsymbol{H}_{B}^B = \boldsymbol{I}_B^B \boldsymbol{\omega} 
    \label{eq1-HBB}
    \end{equation}
        
        (2) The angular momentum of the platform's center of mass relative to the system's overall center of mass;
        \begin{equation}
        \boldsymbol{H}_{BC}^S = \boldsymbol{I}_{BC}^S \boldsymbol{\omega} 
    \label{eq2-HBS}
    \end{equation}
    
        (3) The angular momentum of the rotor with respect to its own center of mass;
         \begin{equation}
        \boldsymbol{H}_{A}^A = \boldsymbol{I}_A^A (\boldsymbol{\omega+\Omega} )
    \label{eq3-HAA}
    \end{equation}
    
        (4) The angular momentum of the rotor's center of mass relative to the system's overall center of mass.
           \begin{equation}
        \boldsymbol{H}_{AC}^S = \boldsymbol{I}_{AC}^S (\boldsymbol{\omega+\Omega} )
    \label{eq4-HAS}
    \end{equation}
     
    Here, $\mathbf{I}_B$ and $\boldsymbol{I}_A$ represent the moments of inertia of the platform and rotor about their mass centers, respectively. Introducing the composite relative moment of inertia $\boldsymbol{I}_S^S=\boldsymbol{I}_{AC}^S+\boldsymbol{I}_{BC}^S$, the equivalant spin inertia tensor $\boldsymbol{I}_1=\boldsymbol{I}_S^S+\boldsymbol{I}_A^A$, The total angular momentum of the spacecraft system about the mass center \(O_S\) can be expressed as Eq. (\ref{eq5-total-angular-momentum}).
    \begin{equation}
        \boldsymbol{H}_{\Sigma}^S = \boldsymbol{I}_B^B \boldsymbol{\omega} +\boldsymbol{I}_1 \left( \boldsymbol{\omega} + \boldsymbol{\Omega} \right)
    \label{eq5-total-angular-momentum}
    \end{equation}
   
   Since the system is free from external torques, Within the inertial reference frame \(OXYZ\), the following are maintained:
    \begin{equation}
        {}^{(N)}\dot{\boldsymbol{H}}_\Sigma^S =  {}^{(B)}\dot{\boldsymbol{H}}_\Sigma^S + \boldsymbol{\omega} \times {}^{(B)}\boldsymbol{H}_\Sigma^S
    \label{eq6-dot-HN}
    \end{equation}
     Combining the above Eqs. (\ref{eq5-total-angular-momentum}) and (\ref{eq6-dot-HN}), the attitude dynamic equation of the system in the body-fixed frame \(O_S X_S Y_S Z_S\) can be derived.
    \begin{equation}
    \begin{split} 
         \dot{\boldsymbol{I}}_1  \left( \boldsymbol\omega + \boldsymbol\Omega \right) + (\boldsymbol{I}_1 +\boldsymbol{I}_B) \dot{\boldsymbol\omega}   + \boldsymbol\omega \times \left[\boldsymbol{I}_1 \left( \boldsymbol\omega + \boldsymbol\Omega \right) \right] + \boldsymbol\omega \times \boldsymbol{I}_B \boldsymbol\omega = 0
    \end{split}        
    \label{eq7-dneq}
    \end{equation}
    
    Herein, the system's total moment of inertia is denoted \( \boldsymbol{I}_2 = \boldsymbol{I}_1 + \boldsymbol{I}_B \). By introducing the dimensionless time \( \tau = \left |  \Omega \right | t \), the influence of the absolute time scale on the system's response is eliminated. Eq.~(\ref{eq7-dneq}) can be reformulated with respect to the dimensionless time scale $\tau$ as:
    
    \begin{equation}
        \dot{\hat{\boldsymbol{\omega}}} = -\boldsymbol{I}_2^{-1} \dot{\boldsymbol{I}}_1 (\hat{\boldsymbol{\omega}} + \boldsymbol{\Omega}_0) - \boldsymbol{I}_2^{-1} \hat{\boldsymbol{\omega}} \times \left[ \boldsymbol{I}_1 (\hat{\boldsymbol{\omega}} + \boldsymbol{\Omega}_0)\right] - \boldsymbol{I}_2^{-1} \hat{\boldsymbol{\omega}} \times \boldsymbol{I}_B \hat{\boldsymbol{\omega}}
    \label{eq8-dot-omega}
    \end{equation}
    Note that \( \hat{\boldsymbol{\omega}} \equiv \frac{\boldsymbol{\omega} }{ \left | \boldsymbol{\Omega}\right |}  \) is the dimensionless angular velocity, defined with respect to the scaled time variable \( \tau = \left | \boldsymbol{\Omega} \right |  t \), and  \( \boldsymbol{\Omega}_0 = \frac{\boldsymbol{\Omega}} { \left | \boldsymbol{\Omega}\right |} \) corresponds to the dimensionless relative angular velocity between the platform and the rotor. From Eq. (\ref{eq8-dot-omega}), it is evident that when  $ {\boldsymbol{\hat{\omega}}}|_{t=0} = 0$, the response \( \hat{\boldsymbol{\omega}} \) of the dynamic equation is fundamentally governed by the inertia coefficients \( \boldsymbol{I}_1 \) and \( \boldsymbol{I}_B \). Variations in \( \boldsymbol{I}_1 \) and \( \boldsymbol{I}_B \) directly influence the characteristics of the dynamic response \( \hat{\boldsymbol{\omega}} \). Consequently, a comprehensive investigation of the impact of inertia parameters \(\boldsymbol{I}_1 \) and \( \boldsymbol{I}_B \) on the dynamic equation's solution provides critical insights for spacecraft design, particularly for systems with large rotors.

\section{Analytical Solutions and Stability Analysis}
\label{sec3-analyticalSol-stability}
    \subsection{Perturbative Reduction and Modal Decoupling of the Dynamic Equations}
    \label{subsec3-1}

To further explore the solution of attitude dynamics and the stability of partial-spin spacecraft, we project Eq.~(\ref{eq8-dot-omega}) onto the spacecraft platform’s body-fixed frame for attitude dynamics computation. In the body-fixed frame of the spacecraft platform, matrix elements of inertia tensors $\boldsymbol{I}_A$, $\boldsymbol{I}_B$ and $\boldsymbol{I}_S$ can be determined as follows:
    \begin{equation}
        \boldsymbol{I}_A = \boldsymbol{T} \cdot \left( \begin{matrix}
        I_{Axx} & 0 & 0 \\
        0 & I_{Ayy} & 0 \\
        0 & 0 & I_{Azz}\\
        \end{matrix} \right)\cdot \boldsymbol{T^{-1}}
    \label{eq9-IA}
    \end{equation}
    \begin{equation}
        \boldsymbol{I}_B = \left( \begin{matrix}
        I_{B_R} & 0 & 0 \\
        0 & I_{B_Y} & 0 \\
        0 & 0 & I_{B_R}\\
        \end{matrix} \right)
    \label{eq10-IB}
    \end{equation}
    \begin{equation}
        \boldsymbol{I}_S = \boldsymbol{T}\cdot\frac{m_B m_A }{(m_A+m_B)}\left( \begin{matrix}
        {h^2} & -dh & 0 \\
        - dh & d^2& 0 \\
        0 & 0 & d^2+h^2\\
        \end{matrix} \right)\cdot \boldsymbol{T^{-1}}
    \label{eq11-IS}
    \end{equation}

    Where $\boldsymbol{T}$ is the rotation matrix between the body-fixed frame of the platform and the rotor, and can be expressed as:
      \begin{equation}
        \boldsymbol{T} =
       \begin{pmatrix}
            cos(\tau) & 0 & sin(\tau) \\
            0 & 1 & 0 \\
            -sin(\tau) & 0 & cos(\tau)
        \end{pmatrix}
        \label{eq12-T}
    \end{equation}

    Thus, the equivalent spin inertia tensor \( \boldsymbol{I}_1 \) in the body-fixed frame of the spacecraft platform can be expressed as:
    \begin{equation}
        \boldsymbol{I}_1 =
        \boldsymbol{T}\cdot\begin{pmatrix}
            I_{xx} & I_{xy} & 0 \\
            I_{xy} & I_{yy} & 0 \\
            0 & 0 & I_{zz}
        \end{pmatrix}\cdot \boldsymbol{T^{-1}}
        \label{eq13-I1}
    \end{equation}

    For simplicity, the matrix coefficients in inertia \( \boldsymbol{I}_1 \) will be used in all subsequent discussions instead of the elements in \( \boldsymbol{I}_A \) and \( \boldsymbol{I}_S \). Introducing a perturbative expansion based on the weak nonlinearity arising from the product of inertia \( I_{xy} \): 
\[
\boldsymbol{\delta} = \mathrm{diag}(\gamma, \gamma^2, \gamma)
\]
where
\[
\gamma = \frac{I_{xy}}{I_{B_R} + I_{zz}}
\]

We then assume the angular velocity admits a perturbative expansion of the form:
\[
    \hat{\boldsymbol{\omega}} = \hat{\boldsymbol{\omega}}^{(0)}+ \boldsymbol{\delta} \cdot \hat{\boldsymbol{\omega}}^{(1)} + \boldsymbol{\delta}^2 \cdot \hat{\boldsymbol{\omega}}^{(2)} + \cdots
\]

Substituting this ansatz into Eq.~(\ref{eq8-dot-omega}) and collecting terms by the order of \( \boldsymbol{\delta} \), we obtain a hierarchy of differential equations. At first order \( \mathcal{O}(\boldsymbol{\delta}) \), the resulting system corresponds to the simplified form presented in Eq.~(\ref{eq9-dot-omega-expanding}). 
    \begin{equation}
    \begin{aligned}
        \dot{{\hat{\boldsymbol{\omega}}}}^{(1)}  =
        \begin{pmatrix}
        \alpha \sin(2\tau) \hat{\omega}_x^{(1)} + \left(\alpha \cos(2\tau) + (\alpha + \beta)\right) \hat{\omega}_z^{(1)} + \sin(\tau) \\
        c_1 \left( 2 \cos(2\tau) \hat{\omega}_x^{(1)} \hat{\omega}_z^{(1)} + \sin(2\tau){(\hat{\omega}_x^{(1)2}} - {\hat{\omega}_x^{(1)2}}) \right) + c_2 \left( \sin(\tau) \hat{\omega}_x^{(1)} + \cos(\tau) \hat{\omega}_z ^{(1)}\right) \\
        \left( \alpha \cos(2\tau) - (\alpha + \beta) \right) \hat{\omega}_x^{(1)} - \alpha \sin(2\tau) \hat{\omega}_z^{(1)} +  \cos(\tau)
        \end{pmatrix}
    \end{aligned}        
    \label{eq9-dot-omega-expanding}
    \end{equation}
    where 
    \begin{equation}  
        \alpha=\frac{\left({I}_{xx}-{I}_{zz}\right)\left({I}_{xx}-{I}_{yy}+2{I}_{B_R}+{I}_{zz}\right)}{2\left({I}_{B_R}+{I}_{zz}\right)\left({I}_{xx}+{I}_{B_R}\right)}
        \label{eq10-alpha}
    \end{equation}
    
    \begin{equation}  
        \beta=-\frac{{I}_{xx}-{I}_{yy}-{I}_{zz}}{{I}_{B_R}+{I}_{zz}}
        \label{eq11-beta}
    \end{equation}

    \begin{equation}  
        c_1 = -\frac{{I}_{xx} - {I}_{zz}}{2\left({I}_{B_Y}+{I}_{yy}\right)}
        \label{eq13-c1}
    \end{equation}
   
    \begin{equation}  
        c_2 = -\frac{(I_{zz} + I_{B_R})(I_{xx} + I_{yy} - I_{zz})}{(I_{B_Y} + I_{yy})(I_{xx} + I_{B_R})}.
        \label{eq18-c2}
    \end{equation}

Eq.~(\ref{eq9-dot-omega-expanding}) reveals that the dynamics in the x- and z-directions are decoupled from that in the y-direction. Accordingly, the attitude dynamics equations along the x- and z-axes can be extracted separately to yield a two-degree-of-freedom system, where $\omega_x$ and $\omega_z$ serve as the state variables:
\begin{equation}
    \left( \begin{matrix}
     \dot{\hat{\omega}}_x^{(1)} \\
    \dot{\hat{\omega}}_z^{(1)} \\
    \end{matrix} \right) = \left(\begin{matrix}
        \alpha\sin(2\tau) &  \alpha\cos(2\tau)+(\alpha+\beta)\\
        \alpha\cos(2\tau)-(\alpha+\beta) & -\alpha\sin(2\tau)\\
        \end{matrix} \right) \left( \begin{matrix}
     \hat{\omega}_x^{(1)} \\
    \hat{\omega}_z^{(1)} \\
    \end{matrix} \right) + \left( \begin{matrix}
    sin(\tau) \\
    cos(\tau) \\
    \end{matrix} \right)\label{eqomegaxy}
\end{equation}

\subsection{Stability analysis of spacecraft attitude motion}
    \begin{theorem}[Floquet theory]\label{thm1}
    Suppose $X \in \mathbb{R}^n$ is a solution to the ordinary differential equation (ODE) of the following LPTV system defined by the T-periodic matrix $\mathbf{A}(t)$. Then, according to Floquet theory, the transition matrix of $\mathbf{\Phi}_{A}(t,t_0)$ can be decomposed into
    \begin{equation}
        \mathbf{\Phi}_{A}(t,t_0)=\mathbf{P}(t) e^{(t-t_0)\mathbf{R}} \mathbf{P}^{-1}(t_0)
    \label{eq15-Phi-A}
    \end{equation}
    where $\mathbf{R}$ is a constant matrix for all $t$ and $\mathbf{P}(t)$ is a non-singular T-periodic matrix.
    \end{theorem}
    For Eq. (\ref{eqomegaxy}), it can be found that $\mathbf{P}(\tau)$ and $\mathbf{R}$ can be expressed as 
    \begin{equation}
        \mathbf{P}(\tau) =
        \left(\begin{matrix}
        \cos \left(\tau\right) & \sin \left(\tau\right)\\
        -\sin \left(\tau\right) & \cos \left(\tau\right)\\
        \end{matrix} \right)
        \label{eq16-Pt}
    \end{equation}
    \begin{equation}
        \mathbf{R} = \left(
        \begin{matrix} 
        0 & u_1 \\
        u_2 & 0 \\
    \end{matrix}\right)
    \label{eq17-R}
    \end{equation}
    where
    \begin{equation}  
        u_1=2\alpha+\beta-1
    \label{eq18-u1}
    \end{equation}
    
    \begin{equation} 
        u_2 = 1 - \beta
    \label{eq19-u2}
    \end{equation}
    
 The eigenvalues of \( R \) are given by:
\[
\lambda_{\text{eig}} = \pm \sqrt{u_1 u_2}
\]
To analyze the stability of the system, we examine the sign of the product \( u_1 u_2 \). By substituting the expressions of \( u_1 \) and \( u_2 \) from Eqs. (18) and (19), we obtain:
\[
u_1 u_2 = -\frac{ (I_{yy} - I_{xx}')(I_{yy} - I_{zz}')}{I_{xx}' I_{zz}'}
\]
where the augmented moments of inertia are defined as \( I_{xx}' = I_{xx} + I_{B_R} \) and \( I_{zz}' = I_{zz} + I_{B_R} \).

As both \( I_{xx}' \) and \( I_{zz}' \) are positive, this motivates the definition of a dimensionless stability criterion:
\[
\sigma = - (I_{yy} - I_{xx}')(I_{yy} - I_{zz}')
\]
The sign of \( \sigma \) determines the sign of \( u_1 u_2 \), and thus governs the qualitative stability behavior of the system.

The corresponding stability regimes are classified as follows:

- When \( \sigma < 0 \), the eigenvalues of \( R \) are purely imaginary, and the system exhibits bounded, periodic motion.

- When \( \sigma = 0 \), the eigenvalues vanish, and the system shows linear divergence and is marginally unstable.

- When \( \sigma > 0 \), the eigenvalues are real. The system is exponentially unstable.

 The relationship between the matrix elements of  \( \mathbf{I}_1 \) and \( \mathbf{I}_B \) and the attitude motion stability is shown through Table \ref{tab2-attitude-stability}. 
\begin{table}[h]
    \centering
    \caption{System attitude stability}
    \label{tab2-attitude-stability}
    \begin{tabular}{c c c c} 
        \hline\hline 
         & $I_{zz'} >  I_{yy}$ & $I_{zz'} =  I_{yy}$ & $I_{zz'} < I_{yy}$ \\
        \hline 
         $I_{xx'} > I_{yy}$ & Stable   & Unstable & Unstable \\    		
         $I_{xx'} = I_{yy}$ & Unstable & Unstable & Unstable \\    		
         $I_{xx'} < I_{yy}$ & Unstable & Unstable & Stable \\
        \hline\hline 
    \end{tabular}
\end{table}
\subsection{Analytical solutions of spacecraft dynamic equations}    
    The solution corresponding to the first-order perturbation equation can be divided into the following cases.
    
    (1) when $\sigma<0$, the state transition matrix for the corresponding homogeneous system of eq.(\ref{eqomegaxy}) is given by
    \begin{equation}
        \mathbf{\Phi} (\tau,0) = \left(
        \begin{matrix}
        \cos(\tau) \cos(\lambda \tau) + \frac{u_2}{\lambda} \sin(\tau) \sin(\lambda \tau) & \frac{u_1}{\lambda} \cos(\tau) \sin( \lambda \tau ) + \sin(\tau) \cos(\lambda \tau) \\
        -\sin(\tau) \cos(\lambda \tau ) + \frac{u_2}{\lambda} \cos(\tau) \sin(\lambda \tau) & -\frac{u_1}{\lambda} \sin(\tau) \sin(\lambda \tau) + \cos(\tau) \cos(\lambda \tau)\\
    \end{matrix}\right)
    \end{equation}
    where
    \begin{equation}
    \begin{aligned}
        \lambda & = \sqrt{-u_1u_2} \\
        & = \sqrt{\frac{-\sigma}{I_{xx'} I_{zz'}}}
    \end{aligned}        
    \end{equation}
    
    Considering $\boldsymbol\omega|_{t=0}=0$, the solution of Eq.~(\ref{eqomegaxy}) is given by
    \begin{equation}
    \left( \begin{matrix}
    \hat{\omega}_x^{(1)} \\
    \hat{\omega}_z^{(1)} \\
    \end{matrix} \right)
    =  \frac{1}{\lambda^2} \left(
    \begin{matrix}
    \lambda\sin(\tau) \sin( \lambda \tau) + \cos{\left( \tau\right)} \cos(\lambda \tau) - \cos(\tau) \\
    \lambda \cos(\tau) \sin(\lambda \tau) - \sin(\tau) \cos(\lambda \tau) + \sin(\tau) \\
    \end{matrix}\right)
    \label{eqxzstablesol}
    \end{equation}

    Substituting Eq. (\ref{eqxzstablesol}) into the Eq. (\ref{eq9-dot-omega-expanding}) and integrating, we obtain
    \begin{equation}
        \hat{\omega}_y^{(1)} = \frac{1}{\lambda^4}c_1u_1 \left( \cos{\left( \lambda \tau\right)} - 1\right)^2 - \frac{1}{\lambda^2}\ c_2 (\cos{\left( \lambda \tau \right)} - 1)
    \end{equation}
    
    Therefore, the solution of the attitude angular velocity can be expressed as
    \begin{equation}
    \left(\begin{matrix}
    \hat{\omega}_x^{(1)} \\
    \hat{\omega}_y^{(1)} \\
    \hat{\omega}_z^{(1)} \\
    \end{matrix}\right) = \frac{1}{\lambda^2} \left( \begin{matrix}
    \lambda \sin(\tau) \sin( \lambda \tau) - u_1 \cos(\tau) \left( \cos(\lambda \tau) - 1 \right) \\
    \frac{1}{\lambda^2} c_1 u_1 \left( \cos \left( \lambda \tau \right) - 1 \right)^2 -\left(\cos ( \lambda \tau ) - 1 \right) \\
    \lambda \cos\left( \tau \right) \sin \left( \lambda \tau \right) + u_1 \sin \left( \tau \right) \left( \cos ( \lambda \tau ) - 1 \right)\\
    \end{matrix}\right)
    \label{eq26-sigma0-omega-xyz}
    \end{equation}
    
    (2) When $\sigma=0$, the state transition matrix is:
    \begin{equation}
        \mathbf{\Phi}( \tau,0) = \left(
        \begin{matrix}
        \cos{\left(\tau\right)} & \sin{\left(\tau\right)} \\
        -\sin{\left(\tau\right)} & \cos{\left(\tau\right)} \\
        \end{matrix}\right)
    \end{equation}
    
    The solution of \( \hat{\omega}_x \) and \( \hat{\omega}_z \) is given by:
    \begin{equation}
    \begin{aligned} 
    \left( \begin{matrix}
    \hat{\omega}_x^{(1) }\\
    \hat{\omega}_z^{(1) } \\
    \end{matrix} \right)    
        &  =  \tau
        \left( \begin{matrix}
        \sin(\tau) \\
        \cos(\tau) \\
    \end{matrix} \right)
    \end{aligned}
    \end{equation}

    Substituting this into the expression for \( \hat{\omega}_y \) in eq.(\ref{eq9-dot-omega-expanding}) and integrating it, we obtain:
    \begin{equation}
    \left( \begin{matrix}
    \hat{\omega}_x^{(1) } \\
    \hat{\omega}_y^{(1) } \\
    \hat{\omega}_z^{(1) } \\
    \end{matrix} \right)
    =  \tau 
    \left( \begin{matrix}
    \sin(\tau)\\
    c_2 \tau\\
    \cos(\tau)\\
    \end{matrix} \right)
    \end{equation}

    This gives the time-dependent evolution of the angular velocities \( \hat{\omega}_x \), \( \hat{\omega}_y \), and \( \hat{\omega}_z \), with the term \( c_2 \) governing the dynamics of \( \hat{\omega}_y \).

    (3) When \( \sigma > 0 \), the state transition matrix is:
    
       \begin{equation}
        \mathbf{\Phi}( \tau, 0) =
        \begin{pmatrix}
        \cos(\tau)\cosh(\lambda' \tau) + \frac{u_2}{\lambda'} \sin(\tau) \sinh(\lambda' \tau) & \frac{u_1}{\lambda'} \cos(\tau) \sinh(\lambda' \tau)+ \sin(\tau) \cosh(\lambda' \tau) \\
        -\sin(\tau) \cosh(\lambda' \tau) + \frac{u_2}{\lambda'} \cos(\tau) \sinh(\lambda' \tau) & -\frac{u_1}{\lambda'}\sin(\tau) \sinh(\lambda' \tau) + \cos(\tau) \cosh(\lambda' \tau)
        \end{pmatrix}
    \end{equation} 
    The result can be expressed as:
    \begin{equation}
    \begin{aligned} 
        \begin{pmatrix}
        \hat{\omega}_x^{(1) } \\
        \hat{\omega}_z^{(1) }
        \end{pmatrix}
        = \frac{1}{\lambda’^2} 
        \begin{pmatrix}
        u_1 \left( \sinh(\lambda' \tau) \sin(\tau) - \cosh^2(\lambda' \tau) + \cosh(\lambda' \tau) \right) \cos(\tau) + \lambda' \sin(\tau) \sinh(\lambda' \tau) \\
        u_1 \left( -\sinh(\lambda' \tau) \sin(\tau) + \cosh^2(\lambda' \tau) - \cosh(\lambda' \tau) \right) \sin(\tau) + \lambda' \cos(\tau) \sinh(\lambda' \tau)
        \end{pmatrix}
    \end{aligned}
    \end{equation}
    where
    \begin{equation}
    \begin{aligned}
        \lambda' & = \sqrt{u_1 u_2} \\
    \end{aligned}        
    \end{equation}
    Thus $\hat{\omega}_y$ can be expressed as
    \begin{equation}
      \hat{\omega}_y^{(1) } =\frac{1}{\lambda'^2}(c_1 cosh^2(\lambda' \tau) +(-2 c_1  + c_2 u_2)cosh(\lambda' \tau)+c_1- c_2u_2  )
    \end{equation}
   
    \label{subsec3-2}
    The above results reveal that the spacecraft's attitude motion stability is fundamentally governed by the inertia characteristics of the rotational components and the platform's principal axes. Specifically, the attitude dynamic response manifests in three distinct regimes: stable motion, instability with time-proportional amplitude growth, and instability characterized by exponential amplitude escalation. The inertia product \( I_{xy} \) influences the solution's amplitude, leaving its form and frequency unaltered.

\section{Precession and Nutation Motion of Spacecraft Attitude}
\label{sec4-precession-nutation}
    To further study the motion characteristics of the spacecraft's attitude dynamics, it is necessary not only to determine the attitude angular velocity of the spacecraft in its coordinate system but also to solve for and analyze its attitude angles in absolute space. Using the $Z$-$Y$-$X$ rotation sequence, the transformation relationship between the Euler angles and under the small-angle assumption, the angular velocity can be expressed as:
    \begin{equation}
        \dot{\boldsymbol{\theta}} =\begin{pmatrix}
           1 & 0 & 0 \\
            0 &1 & \theta_y \\
            0 & 0 &1
        \end{pmatrix}\cdot\boldsymbol{\omega}
        \label{eq36-transitionmatrix}
    \end{equation}

    Nutation and precession are fundamental dynamic behaviors of rigid bodies, particularly in spacecraft mechanics. Precession describes the motion that the spin axis of the rigid body rotates around a fixed axis, referred to as the precession axis. Complementary to precession, nutation represents an additional oscillatory motion of the spacecraft's spin axis, whereby the axis experiences periodic deviations from the primary precession trajectory.

    \begin{figure}[H]
        \centering
        \includegraphics[width=0.5\textwidth]{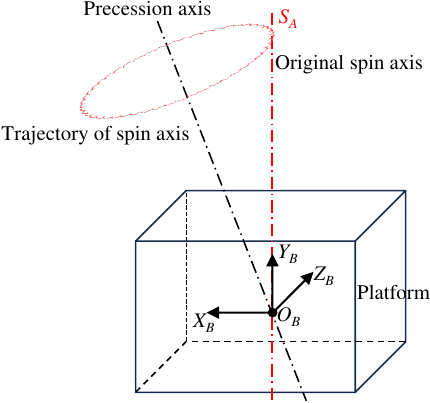}
        \caption{Precession and nutation motion of spacecraft}
        \label{fig3}
    \end{figure}
 \subsection{Nutation and precession motion under stable conditions}
    \label{subsec4-1}
    Based on Eq. (\ref{eq26-sigma0-omega-xyz}) and Eq. (\ref{eq36-transitionmatrix}), the Euler angles of the spacecraft platform \( \theta_x \) and \( \theta_z \) can be derived by integration.
    \begin{equation}
    \left(\begin{matrix}
    \theta_x \\
    \theta_z \\
    \end{matrix}\right)
    =
    \frac{\varepsilon}{\lambda^2-1}
    \left(\begin{matrix}
     \left(-\lambda^2+u_1\right) \sin{\left( \tau\right)}\cos{\left(\lambda \tau\right)}-\lambda\left(u_1-1\right) \cos{\left( \tau\right)} \sin{\left(\lambda \tau\right)} + u_1\left(\lambda^2 - 1\right) \sin{\left( \tau\right)} \\
    \left(-\lambda^2 + u_1\right)\cos{\left( \tau\right)} \cos{\left(\lambda \tau\right)} + \lambda\left(u_1-1\right)\sin{\left( \tau\right)}\sin{\left(\lambda \tau\right)} + u_1\left(\lambda^2 - 1\right) \cos{\left( \tau\right)}-\left(u_1 - 1\right)\lambda^2 \\
    \end{matrix}\right)
    \label{eq37-theta1-3-displacement}
    \end{equation}
   where $\varepsilon=\frac{\gamma}{\lambda^2}$ denotes the amplitude of the first-order solution of the angular velocities $\omega_x$ and $\omega_z$.

    $\theta_y$ can be expressed as
    \begin{equation}
\begin{aligned}
\frac{\gamma}{2\lambda^5  (\lambda^2 - 1)^2} \Bigg[
    & \frac{1}{2} \Bigg(
        \lambda^2 (u_1 - 1)(-\lambda^2 + u_1) \cos(2\tau) 
        - \frac{(\lambda^2 - 1)(u_1 - \lambda)(u_1 + \lambda)}{2}
    \Bigg) \gamma \sin(2\lambda \tau) \\
    & - \frac{\sin(2\tau)(\lambda^4 + (u_1^2 - 4u_1 + 1)\lambda^2 + u_1^2)\lambda\gamma \cos(2\lambda \tau)}{4} \\
    & + (\lambda^2 - 1) \Big[
        u_1 \gamma \lambda^2 (u_1 - 1) \cos(2\tau) 
        + u_1 c_1 \gamma (\lambda^2 - 1) \cos(\lambda \tau) \\
        & - 2 c_2 \lambda^4 
        + (-u_1 (u_1 + 4c_1 + 1)\gamma + 2 c_2)\lambda^2 
        + 2u_1 (u_1 + 2 c_1)\gamma
    \Big] \sin(\lambda \tau) \\
    & + \lambda \Bigg[
        -\frac{\Big(
            (-2u_1\lambda^2 + 2u_1^2)\cos(\lambda \tau) 
            + (u_1^2 + \frac{1}{2})\lambda^2 
            - \frac{3u_1^2}{2}
        \Big)\gamma \sin(2\tau)}{2} \\
        & + \tau \Big(
            2c_2 \lambda^4 
            + \Big(-\frac{1}{2} + u_1^2 + (1 + 3c_1)u_1\Big)\gamma 
            - 2 c_2
        \Big)\lambda^2 
        - \frac{3u_1(u_1 + 2c_1)\gamma}{2}
    \Bigg]
\Bigg]
\end{aligned}
\end{equation}

    From the above results, it can be observed that \( \theta_y \) is smaller by an order of magnitude than \( \theta_x \) and \( \theta_z \), the direction of the relative spin axis $Y_S$ can primarily be described by \( \theta_x \) and \( \theta_z \). Through further analytical manipulation of the Euler angle equations, the relationship between \( \theta_x \) and \( \theta_z \) can be derived:
    \begin{equation}
        \theta_x^2 + (\theta_z + \theta_{z0})^2 = r(\tau)^2 \quad
        \label{eq36-relationship-theta1-3}
    \end{equation}
    where
    \begin{equation}
        \theta_{z0} = \frac{\gamma(u_1 - 1)} {\lambda^2 - 1}
        \label{eq40-theta30}
    \end{equation}
    \begin{equation}
r(\tau) =\sqrt{A_0+A(\tau) }
 \label{eq38-rt}
    \end{equation}
   \begin{equation}
A_0=\frac{\varepsilon^2}{(\lambda^2 - 1)^2}(\frac{ (-\lambda^2 + u_1)^2+\lambda^2(u_1 - 1)^2 }{2} + u_1^2(\lambda^2 - 1)^2)
      \label{eq42}
   \end{equation}
       \begin{equation}
A(\tau)=\frac{\varepsilon^2}{(\lambda^2 - 1)^2} (\frac{(-\lambda^2 + u_1)^2 -\lambda^2(u_1 - 1)^2}{2} \cos(2\lambda \tau) + 2u_1(\lambda^2 - 1)(-\lambda^2 + u_1)\cos(\lambda \tau) )
\label{eq43}
   \end{equation}

    Eq. (\ref{eq36-relationship-theta1-3}) indicates that the attitude motion traces a circle path centered at \( (0, -\theta_{z0}) \), with a dynamically varying radius \( r(\tau) \) that oscillates periodically. This radial motion characterizes the nutation phenomenon, exhibiting two frequencies \( \omega_{n1} = \lambda \) and \( \omega_{n2} = 2\lambda \). 
    The tangential motion of the spin axis on the circle represents precession, with a frequency given by the primary frequency of the Euler angle \( \theta \). $A_0$ denotes the radius of the precession motion, and $A(\tau)$ denotes the change of radius caused by the nutation motion. By the transformation to \( \theta_x \) and \( \theta_z \), the following can be obtained.
    \begin{equation}
        \begin{pmatrix}
        \theta_x \\
        \theta_z
        \end{pmatrix}
        = \frac{\varepsilon}{2(\lambda^2 - 1)}
        \begin{pmatrix}
        (-\lambda^2 + u_1 - \lambda u_1 + \lambda) \sin((1 + \lambda)\tau) + (-\lambda^2 + u_1 + \lambda u_1 - \lambda) \sin((1 - \lambda)\tau) + 2u_1(\lambda^2 - 1) \sin(\tau) \\
        (-\lambda^2 + u_1 + \lambda u_1 - \lambda) \cos((1 + \lambda)\tau) + (-\lambda^2 + u_1 - \lambda u_1 + \lambda) \cos((1 - \lambda)\tau) + 2u_1(\lambda^2 - 1) \cos(\tau) - 2(u_1 - 1)\lambda^2
        \end{pmatrix}
    \end{equation}

    The precession motion can be observed as a superposition of three oscillations with frequencies: \( \omega_{p1} = 1 + \lambda \), \( \omega_{p2} = 1 - \lambda \), and \( \omega_{p3} = 1 \).

    Furthermore, in order to analyze the effect of inertia parameters on nutation behavior during periodic motion, the relative nutation amplitude is defined as
    \begin{equation}
        \varepsilon_n =\sqrt{ \left| \frac{(A(\tau))_{max}}{A_0}\right|} .  
        \label{eq45}
    \end{equation}
   
    The parameter distribution corresponding to different values of \( \varepsilon_n \) in the three-dimensional space formed by \( I_{xx'}, I_{yy}\) , and \(I_{zz'} \) is shown in the figure below:
    
\begin{figure}[H]
    \centering
    \includegraphics[trim = 50mm 100mm 50mm 100mm, clip, width=0.75\textwidth]{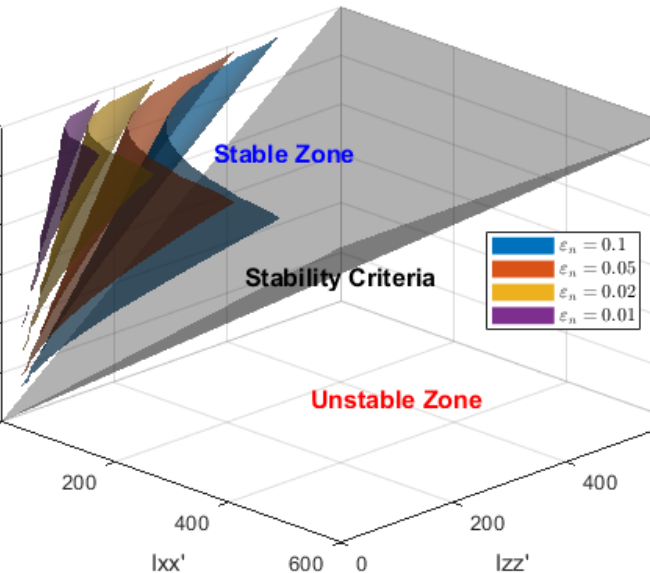}
    \caption{The relationship between relative nutation amplitude $\varepsilon_n$ and inertia parameters}
    \label{fig-peakface}
\end{figure}

    It can be seen that within the common parameter range, as \( I_{xx'} \) and \( I_{zz'} \) decrease and \( I_{yy}\) increases, the relative nutation amplitude decreases. Conversely, as \( I_{xx'} \) and \(I_{zz'} \) increase and \( I_{yy}\) decreases, the relative nutation amplitude increases.
    
\begin{figure}[H]
        \centering
        \includegraphics[width=0.75\textwidth]{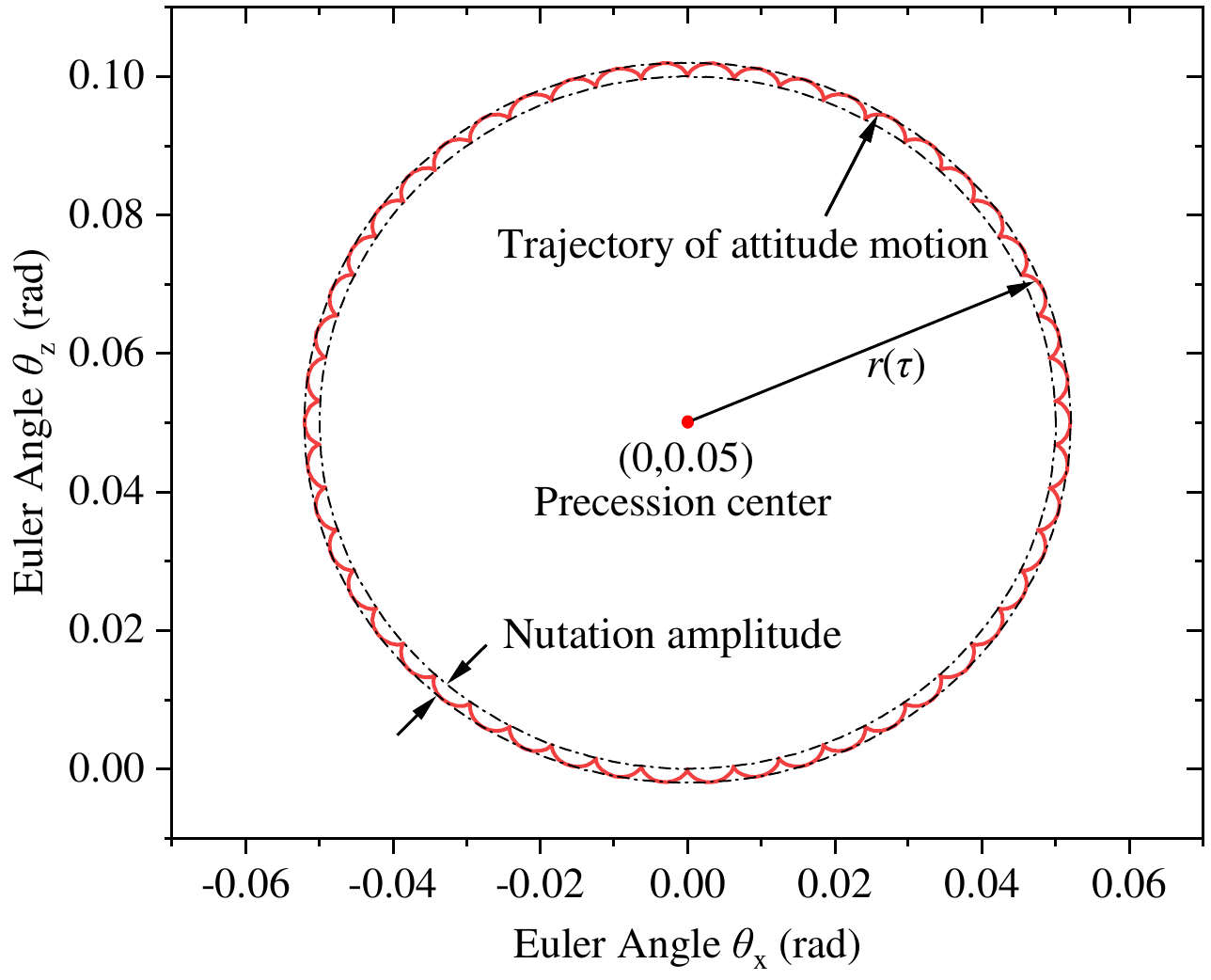}
        \caption{The spacecraft attitude motion on \( \theta_x \) - \( \theta_z \) plane when $\sigma<0$}
        \label{fig-theta1-3}
    \end{figure}

    Fig. \ref{fig-theta1-3} shows the attitude motion of the spacecraft in $\theta_x-\theta_z$ plane over simulation time $T = 320$s in the  $\theta_x - \theta_z$  plane, where the moment of inertia parameters are \( I_{xx} = 1 \) kgm$^2$, \( I_{yy} = 2 \) kgm$^2$, \( I_{zz} = 3 \) kgm$^2$, \( I_{{B_R}} = 100 \) kgm$^2$, and \( I_{xy} = -0.01 \) kgm$^2$, The trajectory reveals a characteristic dynamic behavior wherein the spacecraft's spin axis undergoes periodic precession and nutation.

     \subsection{Nutation and precession motion under unstable conditions}
    \label{subsec4-2}
    When $\sigma=0$, it can be observed from the results above that $\theta_x$ and $\theta_z$ can be expressed as:
    \begin{equation}
    \left(\begin{matrix}
    \theta_x\\
    \theta_z\\
    \end{matrix}\right)
    = \gamma\left(\begin{matrix}
    - \tau + \sin( \tau) \\
    \tau\sin( \tau) + \cos( \tau)-1\\
    \end{matrix}\right)
    \end{equation}
    The relation between $\theta_x$ and $\theta_z$ can be expressed as
    \begin{equation}
        \theta_x^2+\left(\theta_z+\theta_{z0}\right)^2=r\left( \tau\right)^2
    \end{equation}
    where 
    \begin{equation}
        \theta_{z0}=\gamma
    \end{equation}
    \begin{equation}
        r\left( \tau\right)=\left|\gamma\right|\sqrt{\tau^2+1}
    \end{equation}
    This also indicates that the platform's $Y_S$-axis rotates around a fixed center with a radius that varies over time. However, in this case, the precession radius is no longer periodic but grows approximately proportionally with time.

\begin{figure}[htbp]
    \centering
    \begin{subfigure}[b]{0.45\textwidth}
        \includegraphics[width=\textwidth]{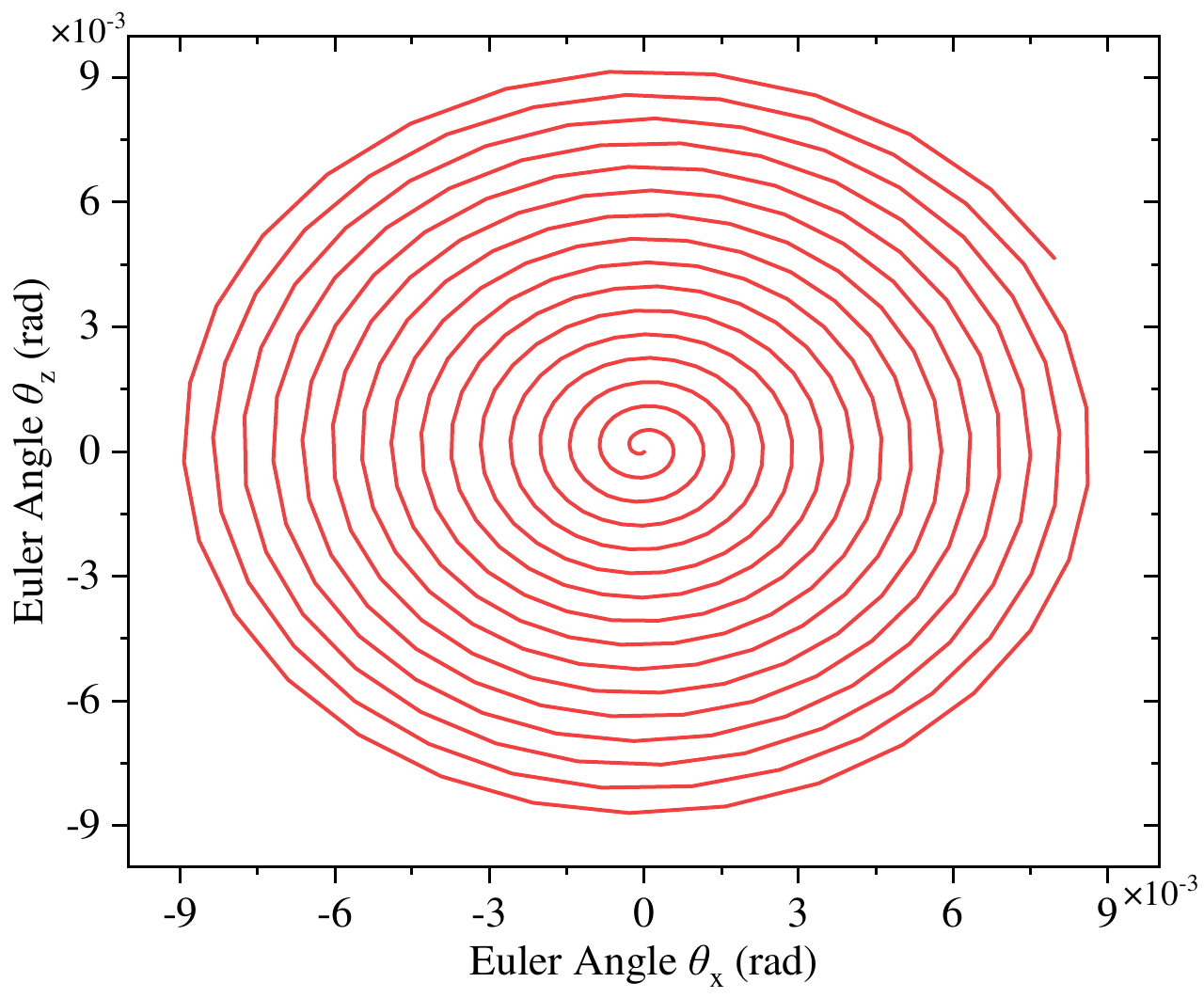}
        \caption{$\sigma=0$}
        \label{fig5-1}
    \end{subfigure}
    \hfill
    \begin{subfigure}[b]{0.45\textwidth}
        \includegraphics[width=\textwidth]{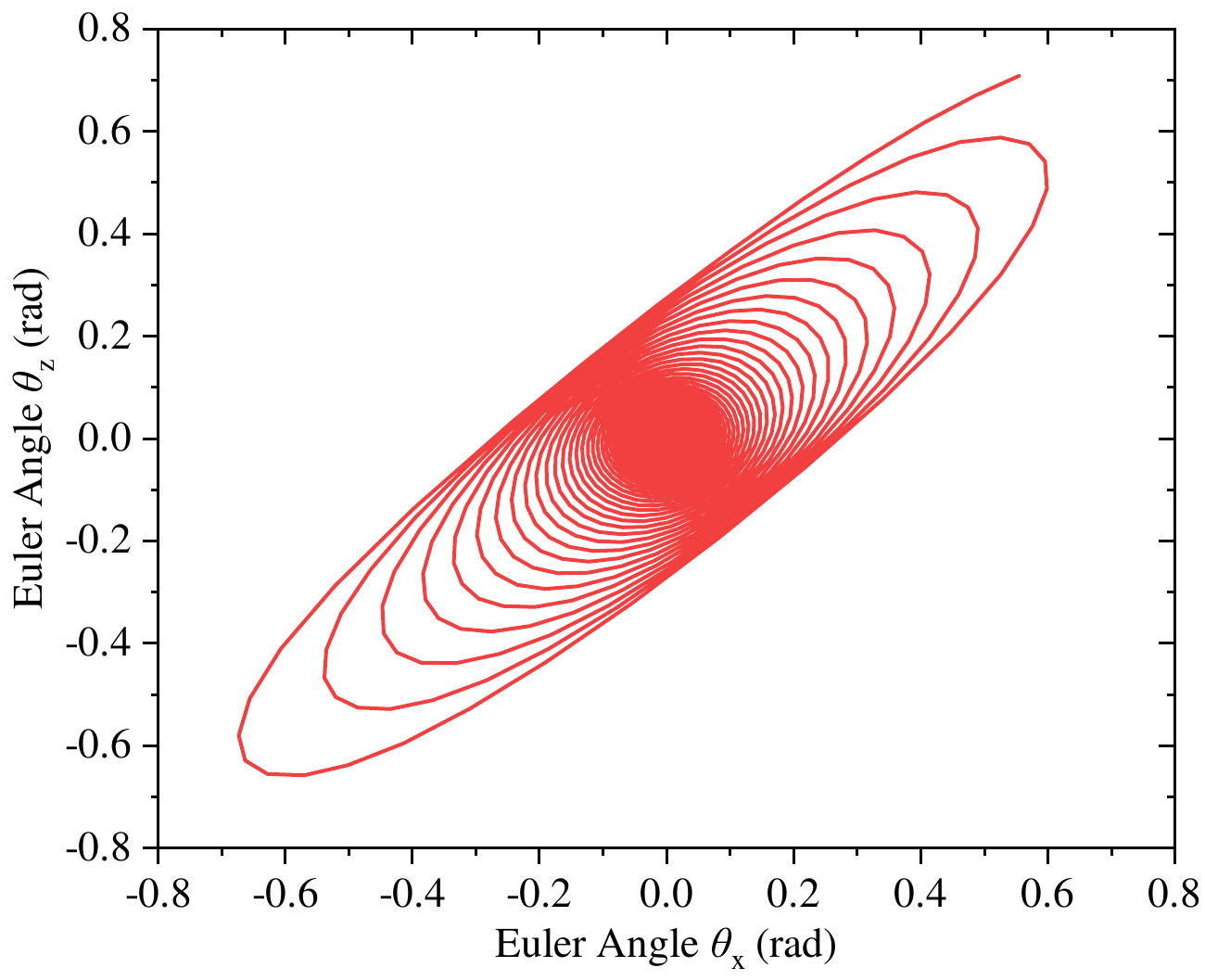}
        \caption{$\sigma>0$}
        \label{fig5-2}
    \end{subfigure}
    \caption{The spacecraft attitude motion on \( \theta_x \) - \( \theta_z \) plane when $\sigma=0$ and $\sigma>0$}
    \label{fig50.}
\end{figure}

     Fig. \ref{fig5-1} shows the projection of the $Y_S$ axis within the platform body frame over simulation time $T = 100s$ in the  $\theta_x - \theta_z$  plane, where the inertia parameters are \( I_{xx} = 3 \) kgm$^2$, \( I_{yy} = 102 \) kgm$^2$, \( I_{zz} = 2 \) kgm$^2$, \( I_{{B_R}} = 100 \) kgm$^2$, and \( I_{xy} = -0.01 \) kgm$^2$. The simulation results validate our conclusions well.

    When \( \sigma > 0 \), attitude angular $\theta_x$ and $\theta_z$ grows exponentially with time. Considering the inertia parameters \( I_{xx} = 1 \) kgm$^2$, \( I_{yy} = 102 \) kgm$^2$, \( I_{zz} = 1 \) kgm$^2$, \( I_{{B_R}} = 100 \) kgm$^2$, and \( I_{xy} = -0.01 \) kgm$^2$, and simulation time $T = 400s$, it can be observed from Fig. \ref{fig5-2} shows that the $Y$-axis orientation of the platform quickly deviates from stability.

\section{Error Analysis}
\label{sec5-error-ayalysis}

    This section focuses on the applicability of the analytical solutions. From the above analytical solutions of the simplified linear system, we know that the magnitude of \( \varepsilon \) determines the convergence rate of the analytical solution. Consequently, it also determines the magnitude of the error in the first-order perturbation solution.
    
    To systematically assess the analytical solution's range of applicability derived from the LPTV system, we compute the error of $\boldsymbol{\hat{\omega}}$ between the first-order analytical solution eq.(\ref{eq26-sigma0-omega-xyz}) and the numerical solution of the full dynamic equation eq.(\ref{eq8-dot-omega}). The comprehensive error analysis is conducted across different configuration parameters for \( \varepsilon < 0.01 \), as detailed in Table \ref{tab1-error}. The error metrics include mean relative error (MRE), mean squared error (MSE), root mean squared error (RMSE), and the coefficient of determination ($R^2$) to provide a multifaceted evaluation of the solution's accuracy.

    \begin{table}[H]
    \centering
    \caption{Error between analytical solution and numerical simulation}
    \label{tab1-error}
    \begin{tabular}{cccccc} 
        \hline\hline
        Simulation time & MRE & MSE & RMSE & $R^2$ \\\hline
        t$<$100s & 0.0011 & 1.6455e-10 &  1.2828e-05 & 0.9999956 \\
        t$<$200s & 0.0044 & 1.8695e-09 & 4.3238e-05 & 0.9999530 \\
        t$<$300s & 0.0079 & 4.8486e-09 & 6.9632e-05 & 0.9998900 \\
       \hline\hline
    \end{tabular}
    \end{table}

\begin{figure}[H]
        \centering        \includegraphics[width=0.6\textwidth]{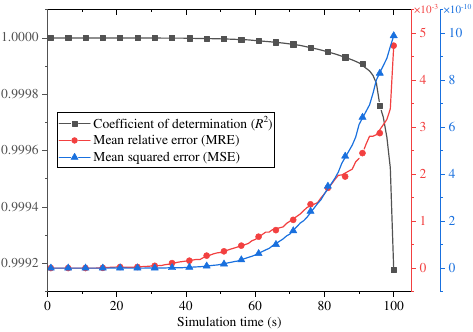}
        \caption{Simulation time vs. error  between analytical solution and numerical simulations}
        \label{fig-error-time}
    \end{figure}

\begin{figure}[H]
        \centering        \includegraphics[width=0.6\textwidth]{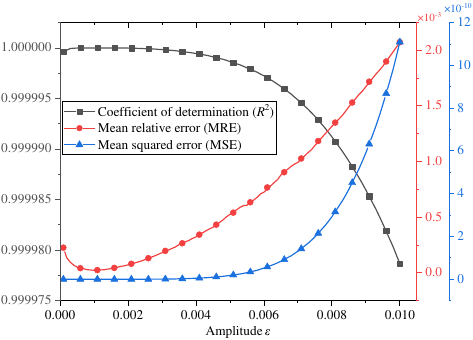}
        \caption{Amplitude of angular velocity $\varepsilon$ vs. error between analytical solution and numerical simulations}
        \label{fig-error-epsilon}
    \end{figure}

    Furthermore, to clarify the relationship between error, simulation duration, and amplitude \( \varepsilon \), time-error and \( \varepsilon \)-error variation curves were calculated under the conditions of simulation duration within 100 s and \( \varepsilon < 0.01 \). These curves are shown in Figs. \ref{fig-error-time} and \ref{fig-error-epsilon}. It can be observed that while the error increases with simulation time and the growth of the parameter \( \varepsilon \), it stays within an acceptable range as long as the preset time and parameter requirements are satisfied. This demonstrates the feasibility of linearizing the dynamic equations under weak system unbalance and the validity of the analytical solution results.

\section{Numerical Examples}
\label{sec6-examples}

\begin{figure}[H]
        \centering
        \includegraphics[width=0.49\textwidth]{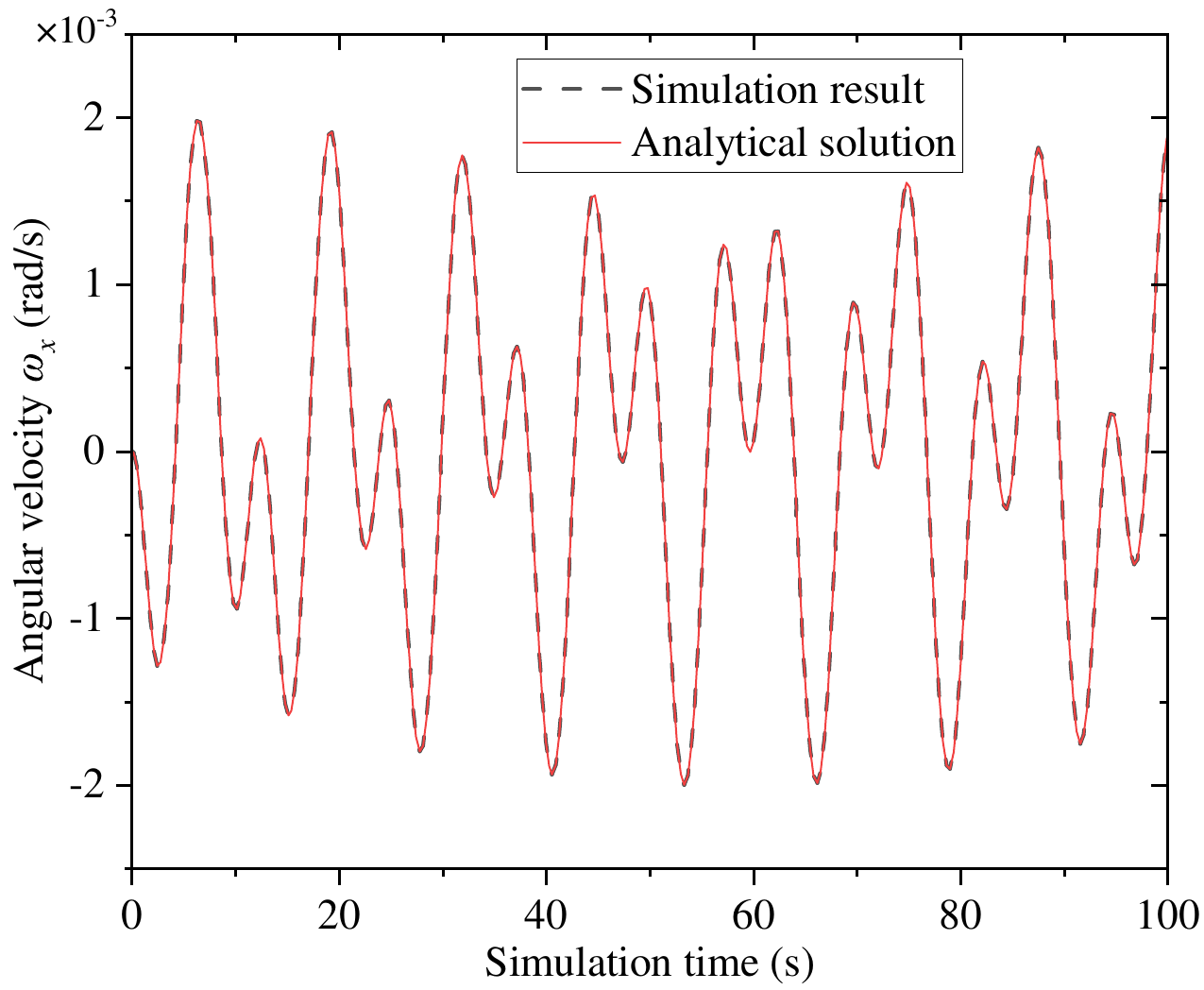}
        \includegraphics[width=0.49\textwidth]{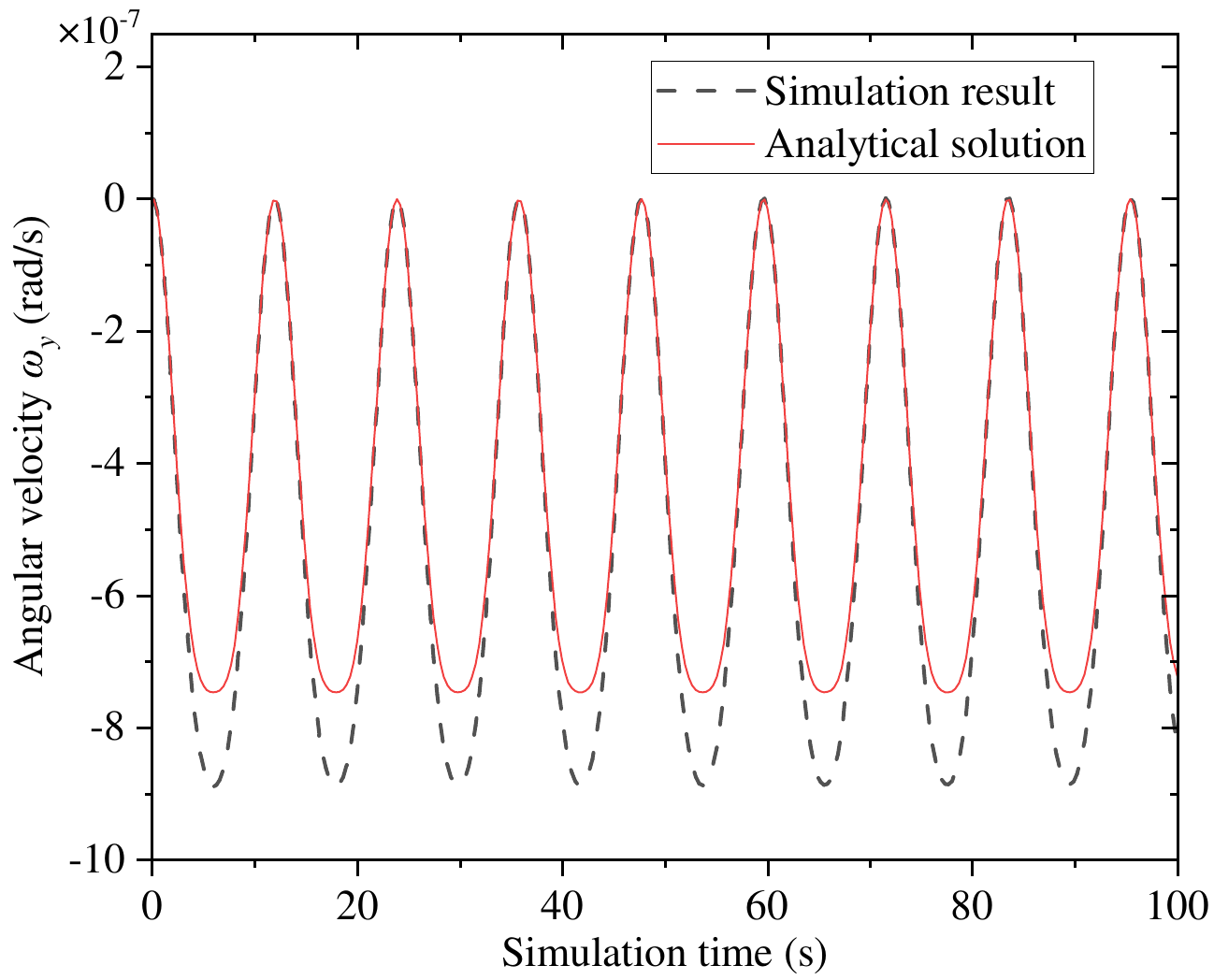}
        \includegraphics[width=0.49\textwidth]{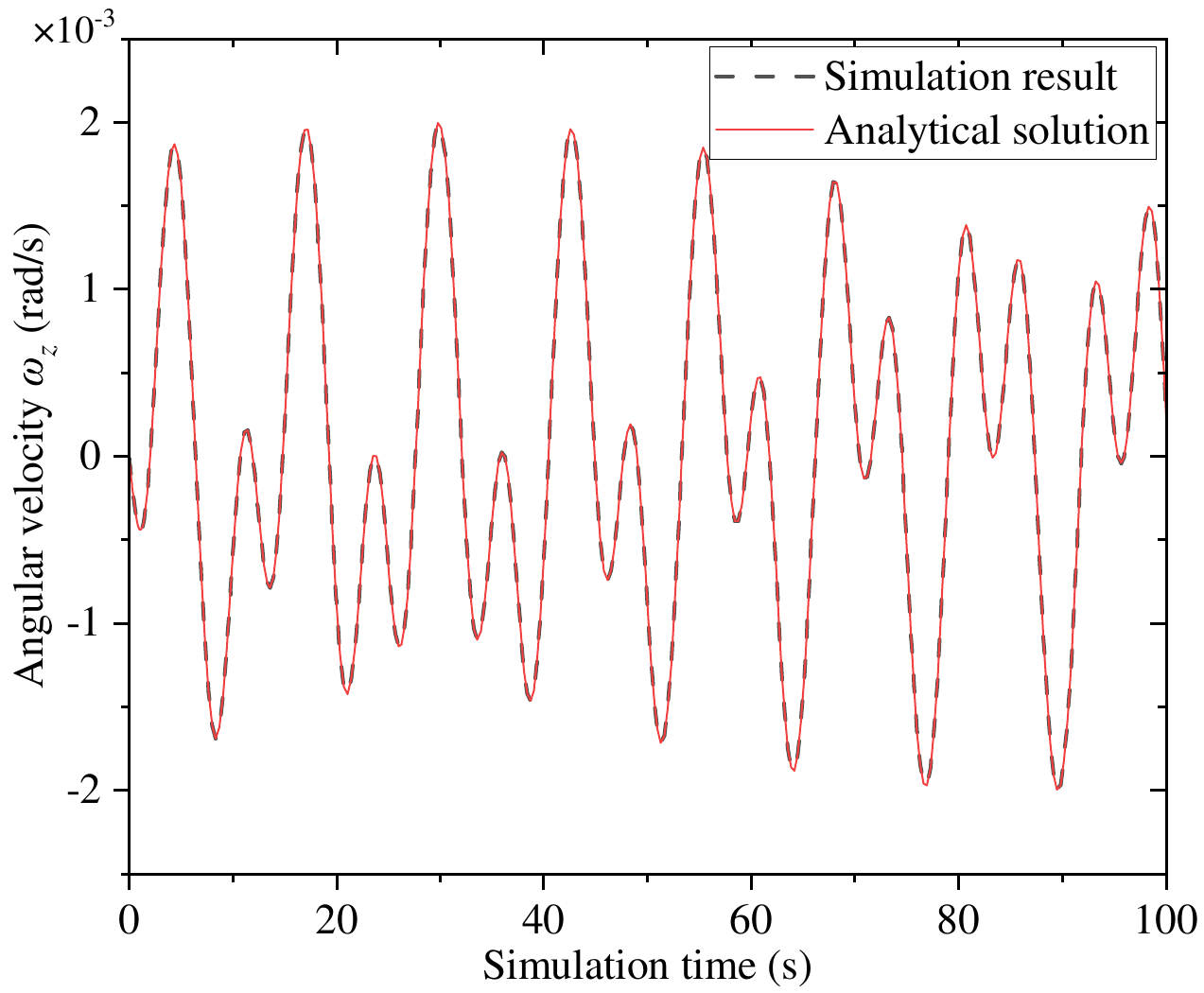}
        \includegraphics[width=0.48\textwidth]{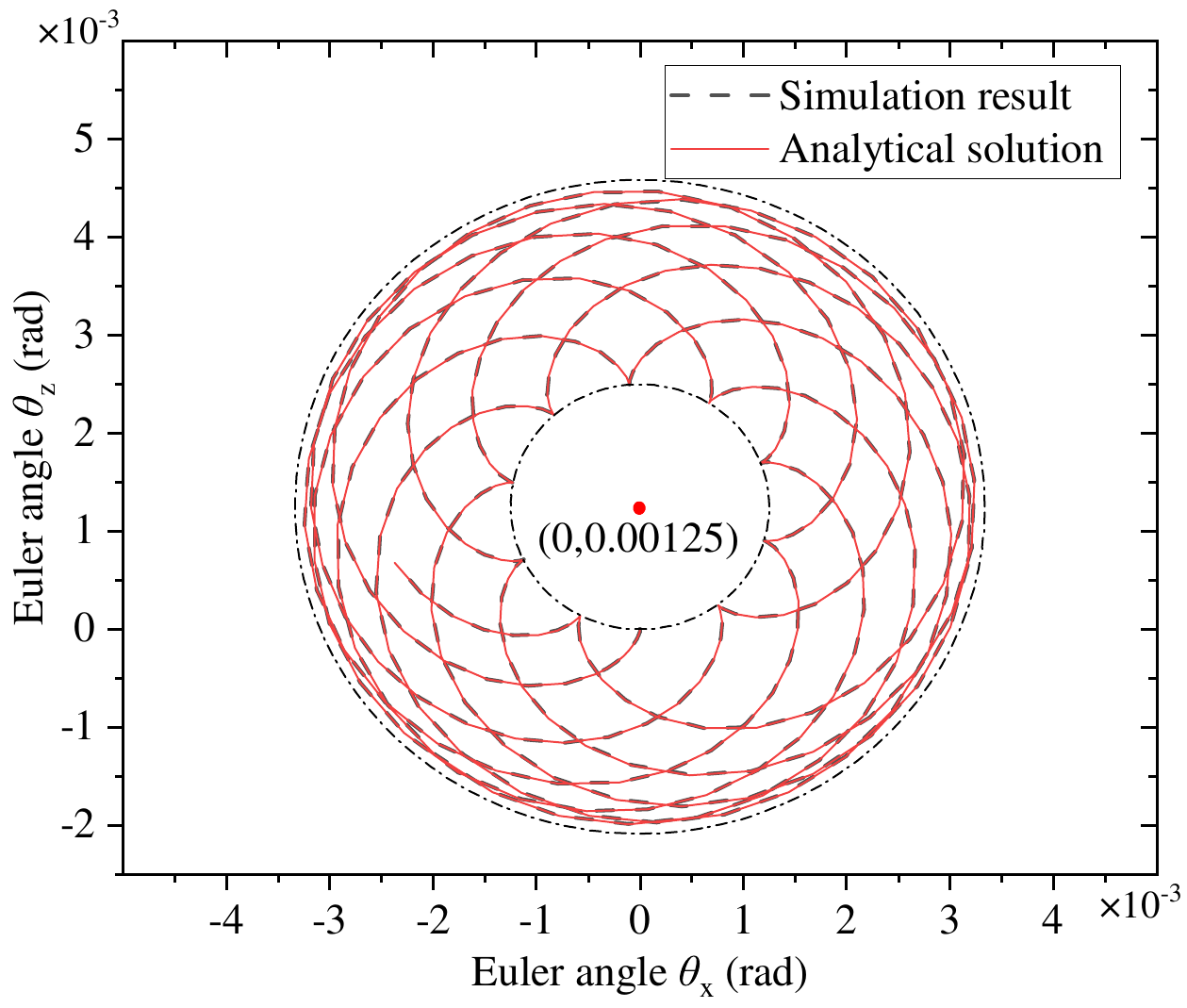}
        \caption{Angular velocity of platform and trajectory of platform orientation ( case 1)}
        \label{fig-exp1-omega-theta13}
    \end{figure}

\begin{figure}[H]
        \centering
        \includegraphics[width=0.49\textwidth]{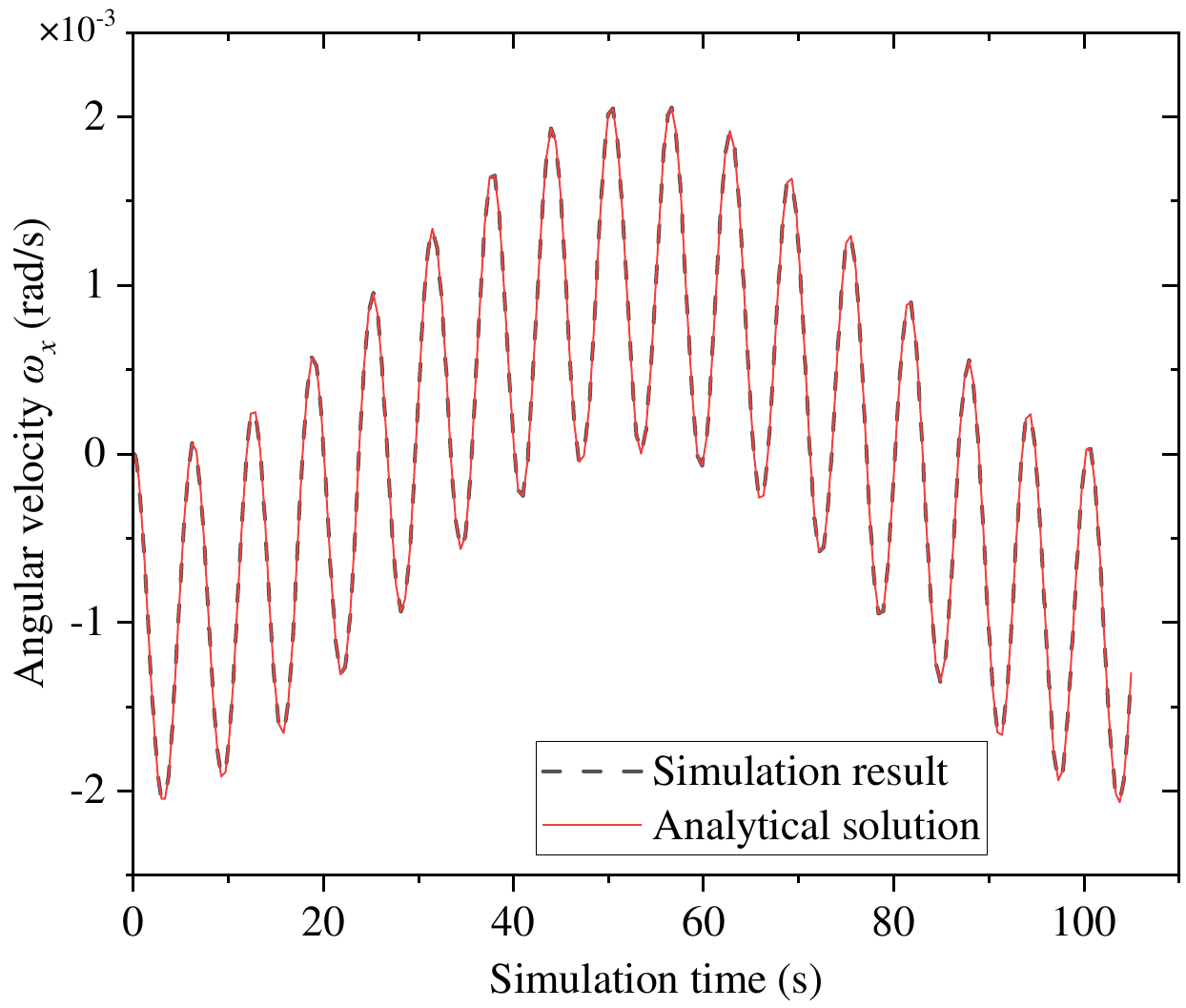}
        \includegraphics[width=0.49\textwidth]{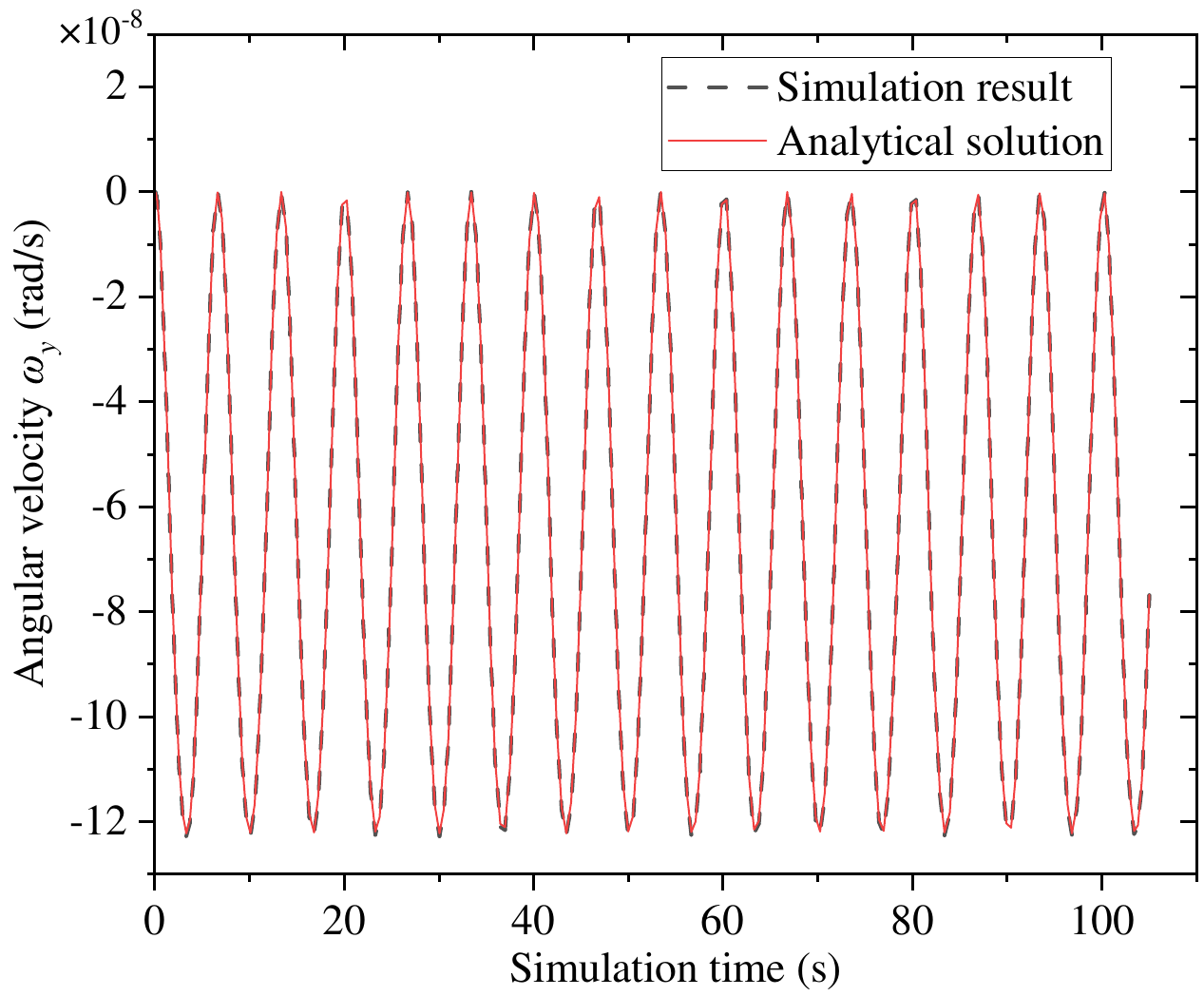}
        \includegraphics[width=0.49\textwidth]{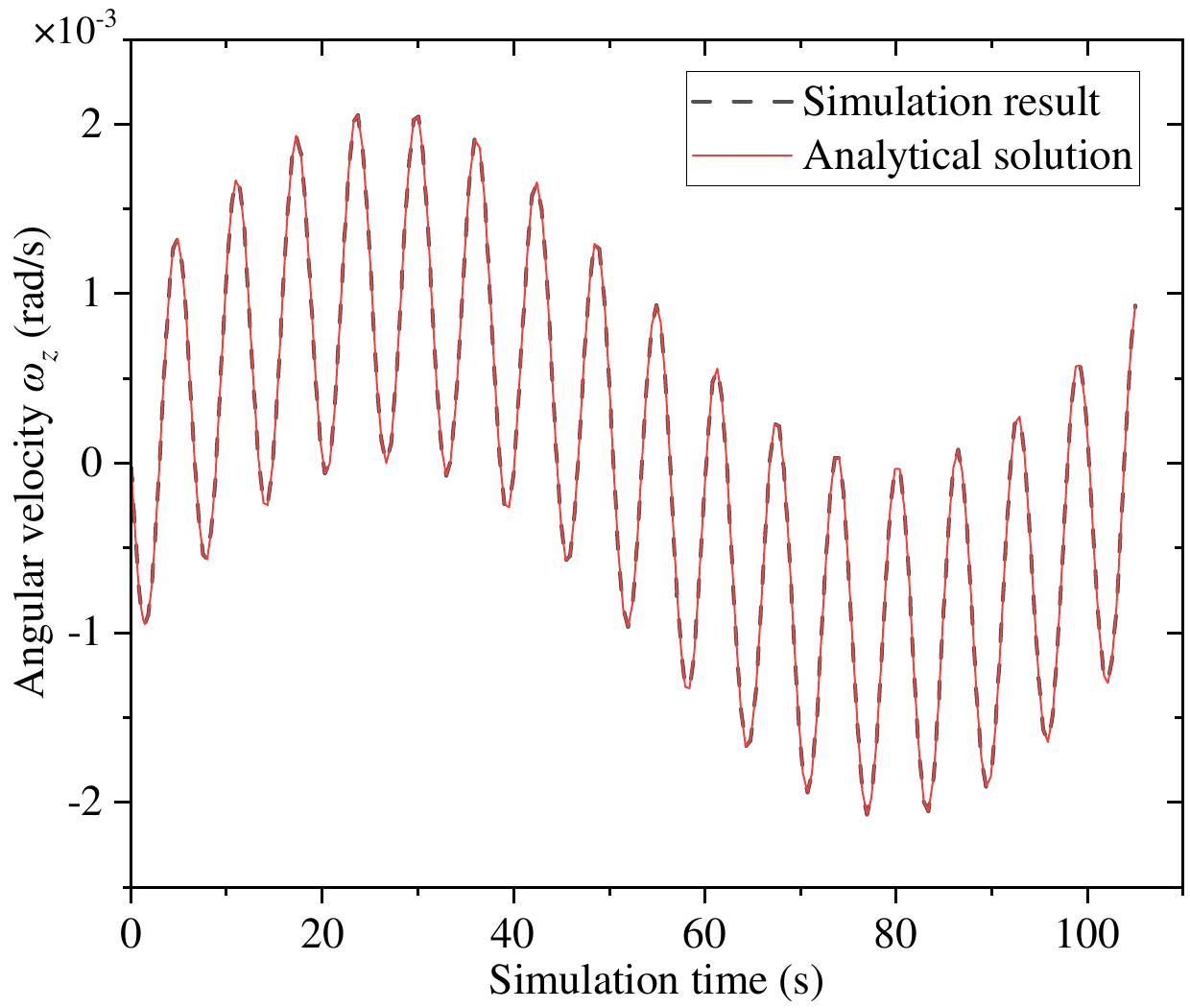}
        \includegraphics[width=0.49\textwidth]{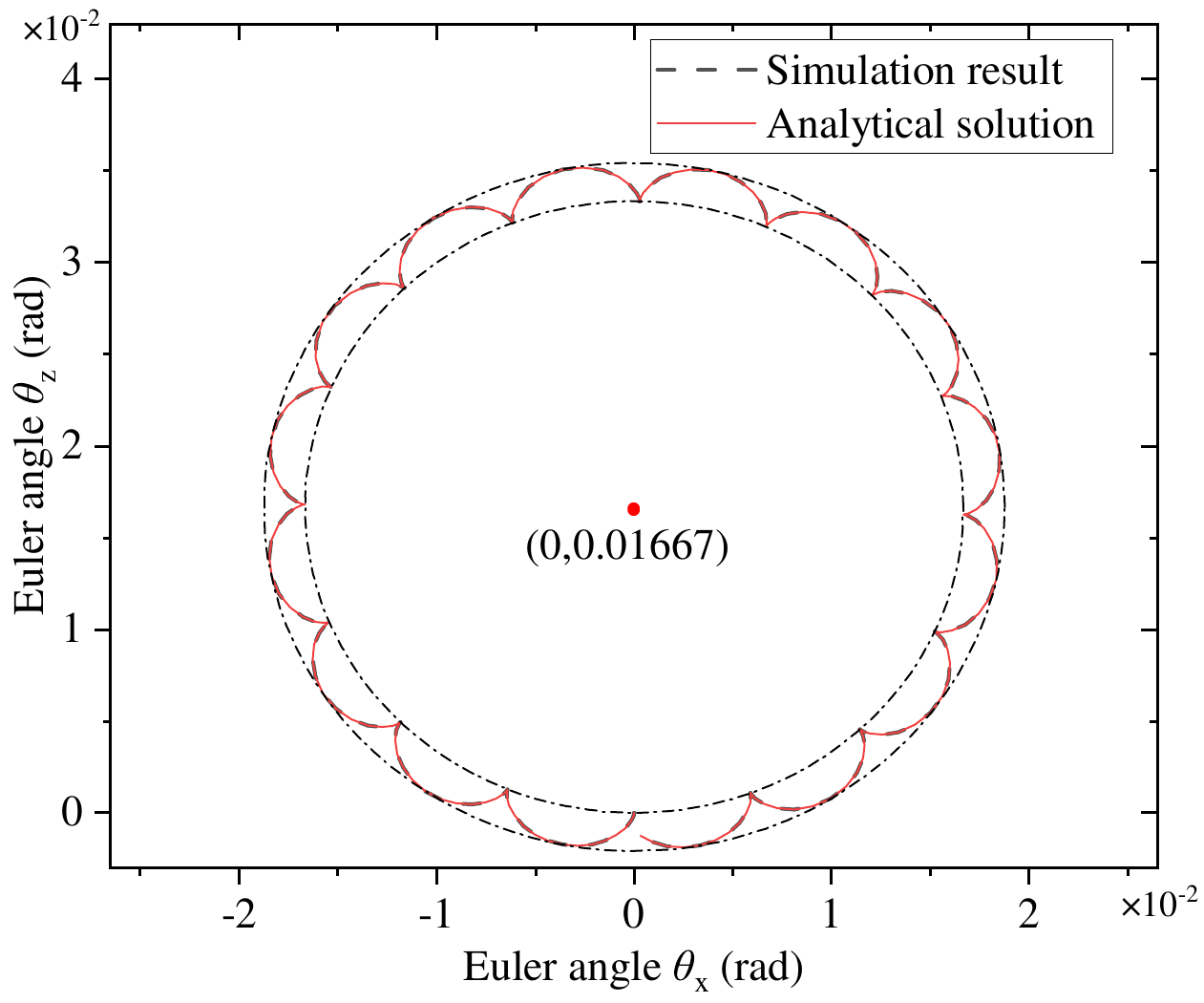}
        \caption{Angular velocity of platform and trajectory of spacecraft platform orientation ( case 2)}
        \label{fig-exp2-omega-theta13}
    \end{figure}

    Herein, we present a comparison between the analytical solutions and the numerical simulations of the system. Two different spacecraft configurations are selected for in-depth analysis. The inertia parameters for each configuration are specified as follows, and the simulation is conducted throughout 100s. The comparative angular velocity curves for the three primary axes and the corresponding trajectory of the platform orientation are illustrated in Figs. \ref{fig-exp1-omega-theta13} and \ref{fig-exp2-omega-theta13}, respectively. The results show that the analytical solutions are in agreement with the original system's simulation results. 

    \begin{itemize}
        \item Example 1: 

    As shown in \cite{feng2022theoretical}, Consider a spacecraft system consisting of a rotor and a platform. Initially, the rotor is dynamically balanced, and its mass center aligns with the principal axis. However, at \( t = 0 \), a portion of the mass at one end of the rotor is separated, disrupting its symmetry. The equivalent rotational moment of inertia can be expressed as: $I_{xx}=80$ kgm$^2$, $I_{yy}=80$ kgm$^2$, $I_{zz}=60$ kgm$^2$, $I_{xy}=-0.1$ kgm$^2$, and the moment of inertia of the platform can be expressed as: $I_{B_R}=100$ kgm$^2$, $I_{B_Y}=90$ kgm$^2$.
 
        \item Example 2: 

        Consider a truss deployed on a satellite platform, rotating relative to the platform about a fixed axis. The length of the truss is 20 m, and its mass is 1 ton. The equivalent rotational moment of inertia can be expressed as: $I_{xx}=20$ kgm$^2$, $I_{yy}=60$ kgm$^2$, $I_{zz}=10$ kgm$^2$, $I_{xy}=-1$ kgm$^2$, and the moment of inertia of the platform can be expressed as: $I_{B_R}=1000$ kgm$^2$,$I_{B_Y}=800$ kgm$^2$.

    \end{itemize}
    
\section{Conclusions}
\label{sec7-conclusion}

    This paper presents a comprehensive analysis of the attitude dynamics of a spacecraft system with an asymmetric platform and an unbalanced rotor. The main achievements are as follows.
    \begin{enumerate}
        \item The dynamic model for asymmetric spacecraft with an unbalanced rotor was developed. The system's nonlinear attitude dynamics were linearized into an LPTV system under weakly unbalanced conditions. Analytical solutions of angular velocities and Euler angles were derived, revealing that the motion frequencies are determined by the primary inertia of the spacecraft platform and the rotor, while the product of the inertia of the rotor relative to the whole system's mass center only affects its amplitude.

        \item A new stability criterion \( \sigma \), composed of the primary inertia of the spacecraft platform and rotor, was established. The sign of \( \sigma \) determines the stability of the attitude motion of the spacecraft system. When \( \sigma < 0 \), the attitude angular velocity and Euler angles exhibit stable and periodic behavior. When \( \sigma \geq 0 \), the system becomes unstable. Specifically, when \( \sigma = 0 \), the attitude angular velocity is proportional to time; and when \( \sigma > 0 \), the attitude angular velocity exponentially diverges over time. This criterion indicates that the stability of spacecraft motion is not affected by the inertia product.

        \item The nutation and precession behavior of the spacecraft dynamics is identified and quantitatively described. Precession occurs around a fixed center on the \( \theta_z \)-axis, while nutation is a periodic oscillation with an amplitude modulated by the spacecraft's inertial properties.

        \item Numerical simulations were conducted to compare the analytical solutions with the numerical solution of the dynamic equation. The results confirmed that the linearized model provides accurate predictions within specific operational domains, particularly when the perturbation parameter $\varepsilon$ is small. The error analysis showed that the linear approximation holds well for $\varepsilon < 0.01$, with the error increasing as $\varepsilon$ grows. This demonstrates the applicability and limitations of the linear model for practical spacecraft design.  
    \end{enumerate}

    The findings of the work are essential for the design and optimization of partial-spin spacecraft. The stability criterion can aid engineers in selecting the appropriate inertia parameters to ensure stable motion. Future work will explore the impact of external torque perturbations and nonlinear couplings on the attitude dynamics, as well as extend the analysis to more complex multi-body spacecraft systems.
    
\section*{Acknowledgments}
    This work is supported by the Key Deployment Project of Chinese Academy of Sciences (22SRM-2023-022) and the Chinese Academy of Sciences Project for Young Scientists in Basic Research (YSBR-107).


\begin{thebibliography}{20}
\newcommand{\enquote}[1]{``#1''}
\providecommand{\natexlab}[1]{#1}
\providecommand{\url}[1]{\texttt{#1}}
\providecommand{\urlprefix}{URL }
\expandafter\ifx\csname urlstyle\endcsname\relax
  \providecommand{\doi}[1]{\discretionary{}{}{}https://doi.org/#1}\else
  \providecommand{\doi}[1]{\discretionary{}{}{}\urlstyle{rm}\url{https://doi.org/#1}}\fi

\bibitem[{Feng et~al.(2022)Feng, Zhang, Zhang, and Li}]{feng2022theoretical}
Feng, G., Zhang, C., Zhang, H., and Li, W., \enquote{Theoretical and Experimental Investigation of Geomagnetic Energy Effect for LEO Debris Deorbiting,} \emph{Aerospace}, Vol.~9, No.~9, 2022, p. 511.
\newblock \doi{10.3390/aerospace9090511}.

\bibitem[{Bracewell and Garriott(1958)}]{bracewell1958rotation}
Bracewell, R.~N., and Garriott, O.~K., \enquote{Rotation of artificial earth satellites,} \emph{Nature}, Vol. 182, No. 4638, 1958, pp. 760--762.
\newblock \doi{10.1038/182760a0}.

\bibitem[{Likins(1966)}]{likins1966effects}
Likins, P.~W., \enquote{Effects of energy dissipation on the free body motions of spacecraft,} Tech. rep., 1966.

\bibitem[{Spencer(1974)}]{spencer1974energy}
Spencer, T., \enquote{Energy-sink analysis for asymmetric dual-spin spacecraft,} \emph{Journal of Spacecraft and Rockets}, Vol.~11, No.~7, 1974, pp. 463--468.
\newblock \doi{10.2514/3.62107}.

\bibitem[{Tsuchiya and Saito(1974)}]{tsuchiya1974dynamics}
Tsuchiya, K., and Saito, H., \enquote{Dynamics of a Spin-Stabilized Satellite Having Flexible Appendages,} \emph{AIAA Journal}, Vol.~12, No.~4, 1974, pp. 490--495.
\newblock \doi{10.2514/3.49274}.

\bibitem[{Likins(2003)}]{likins2003attitude}
Likins, P.~W., \enquote{Attitude stability criteria for dual spin spacecraft,} \emph{Journal of Spacecraft and Rockets}, Vol.~40, No.~6, 2003, pp. 946--951.
\newblock \doi{10.2514/2.7040}.

\bibitem[{Tsuchiya(1979)}]{tsuchiya1979attitude}
Tsuchiya, K., \enquote{Attitude behavior of a dual-spin spacecraft composed of asymmetric bodies,} \emph{Journal of Guidance and Control}, Vol.~2, No.~4, 1979, pp. 328--333.
\newblock \doi{10.2514/3.55883}.

\bibitem[{Cochran et~al.(1982)Cochran, Shu, and Rew}]{cochran1982attitude}
Cochran, J.~E., Shu, P.-H., and Rew, S.~D., \enquote{Attitude motion of asymmetric dual-spin spacecraft,} \emph{Journal of Guidance, Control, and Dynamics}, Vol.~5, No.~1, 1982, pp. 37--42.
\newblock \doi{10.2514/3.56136}.

\bibitem[{Scher and Farrenkopf(1974)}]{scher1974dynamic}
Scher, M., and Farrenkopf, R., \enquote{Dynamic trap states of dual-spin spacecraft,} \emph{AIAA Journal}, Vol.~12, No.~12, 1974, pp. 1721--1725.
\newblock \doi{10.2514/3.49585}.

\bibitem[{Tao and Ramnath(1975)}]{tao1975satellite}
Tao, Y.~C., and Ramnath, R., \enquote{Satellite attitude prediction by multiple time scales method,} Tech. rep., 1975.

\bibitem[{Wang et~al.(2019)Wang, Wang, qin Chen, fei Yue, fei Xie, and peng Chai}]{WANG201991}
Wang, F., Wang, C., qin Chen, X., fei Yue, C., fei Xie, Y., and peng Chai, L., \enquote{High-precision control method for the satellite with large rotating components,} \emph{Aerospace Science and Technology}, Vol.~92, 2019, pp. 91--98.
\newblock \doi{10.1016/j.ast.2019.05.036}.

\bibitem[{Chandra et~al.(2021)Chandra, Kumar, Chattopadhyaya, Chatterjee, and Kumar}]{chandra2021review}
Chandra, M., Kumar, S., Chattopadhyaya, S., Chatterjee, S., and Kumar, P., \enquote{A review on developments of deployable membrane-based reflector antennas,} \emph{Advances in Space Research}, Vol.~68, No.~9, 2021, pp. 3749--3764.
\newblock \doi{10.1016/j.asr.2021.06.051}.

\bibitem[{Schenk et~al.(2014)Schenk, Viquerat, Seffen, and Guest}]{schenk2014review}
Schenk, M., Viquerat, A.~D., Seffen, K.~A., and Guest, S.~D., \enquote{Review of inflatable booms for deployable space structures: packing and rigidization,} \emph{Journal of Spacecraft and Rockets}, Vol.~51, No.~3, 2014, pp. 762--778.
\newblock \doi{10.2514/1.A32598}.

\bibitem[{Bainum et~al.(1970)Bainum, Fuechsel, and Mackison}]{bainum1970motion}
Bainum, P.~M., Fuechsel, P., and Mackison, D., \enquote{Motion and stability of a dual-spin satellite with nutation damping,} \emph{Journal of Spacecraft and Rockets}, Vol.~7, No.~6, 1970, pp. 690--696.
\newblock \doi{10.2514/3.30021}.

\bibitem[{Meng et~al.(2014)Meng, Hao, and Chen}]{meng2014attitude}
Meng, Y., Hao, R., and Chen, Q., \enquote{Attitude stability analysis of a dual-spin spacecraft in halo orbits,} \emph{Acta Astronautica}, Vol.~99, 2014, pp. 318--329.
\newblock \doi{10.1016/j.actaastro.2014.03.001}.

\bibitem[{Janssens and van~der Ha(2011)}]{janssens2011stability}
Janssens, F.~L., and van~der Ha, J.~C., \enquote{On the stability of spinning satellites,} \emph{Acta Astronautica}, Vol.~68, No. 7-8, 2011, pp. 778--789.
\newblock \doi{10.1016/j.actaastro.2010.08.008}.

\bibitem[{Gasbarri et~al.(2016)Gasbarri, Sabatini, and Pisculli}]{gasbarri2016dynamic}
Gasbarri, P., Sabatini, M., and Pisculli, A., \enquote{Dynamic modelling and stability parametric analysis of a flexible spacecraft with fuel slosh,} \emph{Acta Astronautica}, Vol. 127, 2016, pp. 141--159.
\newblock \doi{10.1016/j.actaastro.2016.05.018}.

\bibitem[{Liu et~al.(2018)Liu, Yue, and Zhao}]{liu2018attitude}
Liu, F., Yue, B., and Zhao, L., \enquote{Attitude dynamics and control of spacecraft with a partially filled liquid tank and flexible panels,} \emph{Acta Astronautica}, Vol. 143, 2018, pp. 327--336.
\newblock \doi{10.1016/j.actaastro.2017.11.036}.

\bibitem[{Liu and Chen(2013)}]{liu2013chaos}
Liu, Y., and Chen, L., \emph{Chaos in attitude dynamics of spacecraft}, Springer, 2013.
\newblock \doi{10.1007/978-3-642-30080-6}.


\end{thebibliography}

\end{document}